\RequirePackage{fix-cm}
\documentclass[twocolumn]{svjour3}          
\smartqed  
\usepackage{graphicx}

\usepackage{amsmath}
\usepackage{amssymb}
\usepackage{amsbsy}
\usepackage{float}

\usepackage{subfigure}

\newcommand{\tr}[1]{\mathop{\mathrm{trace}}\left(#1\right)}

\newcommand{\V}[1]{\ensuremath{\boldsymbol{#1}}}
\newcommand{\M}[1]{\ensuremath{\mathsf{#1}}}
\newcommand{\VT}[1]{\ensuremath{\boldsymbol{#1}}^{\textrm{T}}}
\newcommand{\MT}[1]{\ensuremath{\mathsf{#1}}^{\textrm{T}}}
\newcommand{\kth}[1]{\ensuremath{{#1}^{\textrm{th}}}}

\newcommand{\pdiff}[2]{ \frac{\partial #1}{\partial#2} }

\begin{document}

\title{Regularized Reconstruction of a Surface from its Measured Gradient Field}
\subtitle{Algorithms for Spectral, Tikhonov, Constrained, and
Weighted Regularization}


\author{Matthew Harker         \and
        Paul O'Leary 
}


\institute{M. Harker and P. O'Leary \at
              Institute for Automation \\
              University of Leoben \\
              Peter-Tunner-Strasse 27 \\
              8700 Leoben, Austria \\
              Tel.: +43-3842-402 5309\\
              Fax: +43-3842-402 5302\\
              \email{matthew.harker@unileoben.ac.at}\\           
              \email{automation@unileoben.ac.at}           
}

\date{Received: date / Accepted: date}

\maketitle

\begin{abstract}
This paper presents several new algorithms for the regularized
reconstruction of a surface from its measured gradient field.  By
taking a matrix-algebraic approach, we establish general framework
for the regularized reconstruction problem based on the Sylvester
Matrix Equation.  Specifically, Spectral Regularization via
Generalized Fourier Series (e.g., Discrete Cosine Functions, Gram
Polynomials, Haar Functions, etc.), Tikhonov Regularization,
Constrained Regularization by imposing boundary conditions, and
regularization via Weighted Least Squares can all be solved
expediently in the context of the Sylvester Equation framework.
State-of-the-art solutions to this problem are based on sparse
matrix methods, which are no better than
$\mathcal{O}\!\left(n^6\right)$ algorithms for an $m\times n$
surface.  In contrast, the newly proposed methods are based on the
global least squares cost function and are all
$\mathcal{O}\!\left(n^3\right)$ algorithms. In fact, the new
algorithms have the same computational complexity as an SVD of the
same size. The new algorithms are several orders of magnitude
faster than the state-of-the-art; we therefore present, for the
first time, Monte-Carlo simulations demonstrating the statistical
behaviour of the algorithms when subject to various forms of
noise.  We establish methods that yield the lower bound of their
respective cost functions, and therefore represent the
``Gold-Standard" benchmark solutions for the various forms of
noise. The new methods are the first algorithms for regularized
reconstruction on the order of megapixels, which is essential to
methods such as Photometric Stereo. \\
\end{abstract}

\section{Introduction}
%
Surface reconstruction from a gradient field is an important
problem, not only in Imaging, but in the Physical Sciences in
general; it is essential to many applications such as Photometric
Stereo~\cite{woodham1980}, Seismic Imaging~\cite{robein}, as well
as the more general problem of the numerical solution of Partial
Differential Equations.  The reconstruction from gradients problem
can be considered to be an inverse problem, that is, inversion of
the process of differentiation. The difficulty arises in the fact
that if a gradient field is corrupted by Gaussian noise, it is
generally no longer integrable. To make matters worse, the
Gaussian noise is itself to some degree integrable, which
introduces bias into the solution. It is generally known that the
surface can be reconstructed up to a constant of integration,
whereby a global least squares solution accomplishes
this~\cite{Harker2008c}. However, when different forms of noise
are present (e.g., lighting variations in Photometric Stereo, or
gross outliers), the least squares solution is no longer optimal
in the maximum likelihood sense.  To suppress such varied types of
noise, some form of regularization is required on the solution to
the reconstruction problem.  More importantly, a mathematically
sound and efficient solution to this problem is fundamental to
obtaining useable results from surface measurement via Photometric
Stereo. In this paper, we derive several new methods which
incorporate state-of-the-art regularization techniques into the
surface reconstruction problem.
In~\cite{Harker2008c}, it was first shown that the global least
squares minimizer to the reconstruction problem satisfies a
Sylvester Equation, that is, the reconstructed surface $\M{Z}$
satisfies a matrix equation of the form
\begin{equation} \label{eqn:sylEq0}
\M{A}\M{Z} + \M{Z}\M{B} = \M{C}.
\end{equation}
In this paper, it is demonstrated that the Sylvester Equation is
fundamental to the surface reconstruction from gradients problem
in general, in that it leads to $\mathcal{O}\!\left(n^3\right)$
algorithms for all of the most effective forms of regularization
(See, e.g., Engl~\cite{Engl}) for the reconstruction problem.  No
existing method provides a regularized solution with this order of
efficiency. Some preliminary portions of the material herein
appeared in~\cite{Harker2008c,harker2011}.
\subsection{Previous Methods}
The most basic of surface reconstruction algorithms are based on
the fact that the line integral over a closed path on a continuous
surface should be zero.  The line-integral
methods~\cite{wu1988,klette,robleskelly2005}, optimize local
least-squares cost functions, and vary mainly only in the
selection of integration paths, ranging from simple~\cite{wu1988}
to elegant~\cite{robleskelly2005}.  This local nature means that
the reconstruction is only optimal locally over each integration
path. The error residual is non-uniform over the surface, and
hence the methods are not optimal in the presence of Gaussian
noise.  They have the further disadvantage that no form of
regularization (global or otherwise) can be incorporated into the
solution. \\
Horn and Brooks~\cite{horn1986} proposed to take a global approach
to the optimization problem, by means of the Calculus of
Variations.  The problem is formulated in the continuous domain as
the minimization over the domain $\mathcal{D}$ of the functional,
\begin{equation} \label{eqn:varFunc}
J = \iint_{\mathcal{D}} \! \left( z_x - \hat{z}_x \right)^2 +
\left( z_y - \hat{z}_y \right)^2 \, \mathrm{d}x \, \mathrm{d}y
\end{equation}
subject to the boundary conditions,
\begin{equation} \label{eqn:varBCs}
\phi( z, z_x, z_y, x, y, t) = 0,
\end{equation}
whereby $\hat{z}_x$ and $\hat{z}_y$ is the measured gradient. The
solution satisfies the associated Euler-Lagrange equation,
\begin{equation}
\pdiff{^2z}{x^2} + \pdiff{^2z}{y^2} = \pdiff{}{x}\hat{z}_x +
\pdiff{}{y}\hat{z}_y
\end{equation}
which is known as Poisson's Equation.  It should be stressed, that
this equation alone does not specify the solution uniquely; a
unique solution to this boundary value problem is only obtained
when the boundary conditions (or another constraint) are
specified.  They developed an iterative averaging scheme in the
discrete domain with the aim of solving this problem, but they
found it to be non-convergent; however, plausible (but still
biased) results are obtained after several thousand
iterations~\cite{durou2007}.
\\
Further methods were developed based on the variational approach
in the context of shape from shading.  Frankot and
Chellappa~\cite{frankot1988} solved the reconstruction problem of
Equations (\ref{eqn:varFunc}) and (\ref{eqn:varBCs}) by a discrete
Fourier Transform method, whereas Simchony et
al.~\cite{simchony1990} used a Discrete Cosine Transform method
for solving the Poisson Equation~\cite{dorr1970}. The solution of
Frankot and Chellappa assumes periodic boundary
conditions\footnote{Periodic boundary conditions for the
rectangular domain $x \in \left[a,b\right]$, $y \in
\left[c,d\right]$, imply that the surface satisfies
$z(a,y)=z(b,y)$ and $z(x,c)=z(x,d)$.  That is, the surface takes
on the same values at opposing boundaries.  Hence, periodic
boundary conditions are largely unrealistic in real-world
applications such as Photometric Stereo.}, whereby the method
projects the gradient onto complex Fourier Basis functions; the
reconstruction can be accomplished by means of the Fast Fourier
Transform (FFT). The approach of Simchony et al. uses cosine
functions under the assumption that they satisfy homogeneous
Neumann boundary conditions; implementations of the algorithm
unfortunately require zero-padding the gradient and thus introduce
unnecessary bias into the solution. Other methods function in a
similar manner, i.e., by projecting the measured gradient field
onto a set of integrable basis functions. Kovesi~\cite{kovesi2005}
also assumes periodic boundary conditions, but uses shapelets for
basis functions.  In practice, periodic boundary conditions are
unrealistic since, for example, a function of the form $z(x,y) =
ax + by$ is impossible to reconstruct; their results are therefore
mainly only of theoretical interest.  Kara\c{c}al\i\ and Snyder's
method~\cite{karacali2003,karacali2004} effectively uses Dirac
delta functions, but requires the storage and orthogonalization of
a $2mn\times mn$ matrix, and is computationally cumbersome at
best. It should be noted that while basis functions can be used to
solve the integration problem, they have yet to be used for the
purpose of regularization of the
surface reconstruction problem. \\
  Finally, Harker and O'Leary~\cite{Harker2008c} showed that an unconstrained
solution analogous to the integration of a gradient field could be
obtained by working directly in the discrete domain. The global
least squares cost function was  formulated in terms of matrix
algebra; it was shown that the minimizing solution satisfied a
matrix Lyapunov (Sylvester) Equation, and that the solution was
unique up to a constant of integration. This approach represents
the basic least-squares solution on which regularized
least-squares solutions can be based; as such, throughout this
paper it will be referred to as the GLS (global least-squares)
solution.  It is the ``Gold-Standard" benchmark solution when the
gradient field is corrupted by i.i.d.\ Gaussian noise.  As their
methodology is fundamental to the methods derived in this paper,
it is described in more detail in Section \ref{sec:GLSSR}. \\
  As for the state-of-the-art reconstruction methods which incorporate
some form of
regularization~\cite{horovitz2004,agrawal2006,ng2010}, they are
similarly all based on Poisson's Equation, and therefore also
require the specification of some form of boundary conditions.
Their greatest disadvantages are they formulate the optimization
problem by ``vectorizing"~\cite{vanLoan2000} the surface $\M{Z}$,
which involves stacking the $z_{ij}$ into a vector, resulting in a
$2mn\times mn$ coefficient matrix.  Generally, this results in an
$\mathcal{O}\!\left(n^6\right)$ algorithm to solve the linear
system.  Due to their sheer size, sparse iterative methods, such
as LSQR~\cite{paige1982}, must be used.  Those methods with
nonlinear optimization problems~\cite{agrawal2006,ng2010} thus use
iterative methods nested within iterative methods and become
computationally unfeasible with increasing surface size.
With regards to the $\mathcal{O}\!\left(n^6\right)$ algorithms in
general (i.e.,
~\cite{lee1993,karacali2003,horovitz2004,agrawal2006,ng2010,koskulics2012}),
an indication of their impracticality can be gleaned from the
published statistics; some have gargantuan memory requirements and
are limited to surfaces of $32\times 32$~\cite{karacali2003};
others are excruciatingly slow, requiring $3.5$ hours to
reconstruct a $240\times 314$ surface~\cite{ng2010}.  Clearly none
of these methods can be used for any practical purposes, such as
Industrial Photometric Stereo.
\\
From this body of literature aimed at the reconstruction of a
surface from its discrete gradient, we can summarize the following
problems which are, until now, still open problems:
\begin{itemize}
\item Each method solves only one particular sub-problem, e.g.,
reconstruction with boundary conditions.
\item  Moreover, the methods which solve the more complicated
problems, such as regularized reconstruction, are grossly
impractical.
\item Most importantly, each method lack generality, e.g., the
Frankot-Chellappa or Simchony et al. methods do not solve the
Tikhonov Regularization problem.
\end{itemize}
In this paper, we propose a computational framework based on the
Sylvester Equation, which solves all the main regularization
problems which can be associated with the reconstruction of a
surface from its discrete gradient field.
Moreover, all algorithms presented in this paper are shown to be
of $\mathcal{O}\!\left(n^3\right)$ complexity.  To comprehend this
improvement, recall that the development of the FFT reduced a
computation of $\mathcal{O}\!\left(n^2\right)$ to an
$\mathcal{O}\!\left(n\log{n}\right)$ complexity, which for the
computers of the time meant the reduction of a near-impossible
computation to a reasonably efficient computation. The algorithms
presented here represent the first practically applicable
algorithms for the regularized least-squares reconstruction
problem\footnote{The MATLAB code implementing the methods
presented in this paper is available at
http://www.mathworks.com/matlabcentral/fileexchange/
authors/321598}.
%
\section{Global Least Squares Surface Reconstruction}
\label{sec:GLSSR}

\subsection{Numerical Differentiation}

\begin{figure}[!t]
\centering
\includegraphics{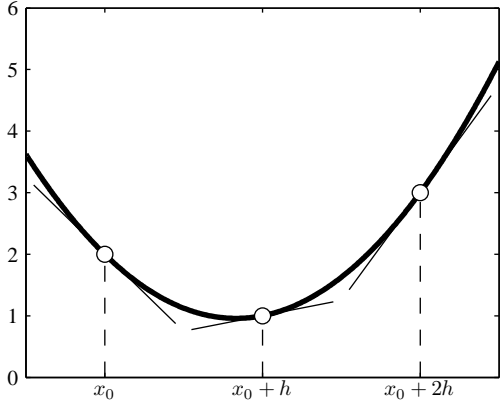}
\caption{Numerical derivatives of a three point sequence.  For the
three point sequence there is a single interpolating parabola. Its
tangent lines at the discrete points are shown to indicate the
numerical derivatives.} \label{fig:threePt}
\end{figure}

\begin{figure}[!t]
\centering
\includegraphics{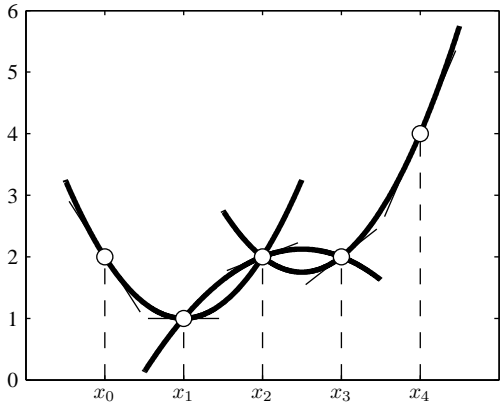}
\caption{Numerical derivatives of a five point sequence with
second order accurate formulas.  There are three interpolating
parabolas, whereby the first and last determine the derivatives
for the end points.  The indicated tangent slopes are all second
order accurate.} \label{fig:fivePt}
\end{figure}

The numerical differentiation of a discrete signal is most
commonly computed by differentiating the polynomials which
interpolate it locally. Using the Lagrange interpolation
polynomials and their corresponding error terms (Lagrange
remainder), one obtains differentiation formulas along with their
error estimates; see Burden and Faires~\cite{Burden}.  For
example, for the three point sequence $\{x_0,x_0+h,x_0+2h\}$ with
even spacing $h$, the derivative of $y=f(x)$ at the first point
$x_0$ is given as,
\begin{eqnarray} \label{eqn:derFormL}
f^{\prime} \! \left(x_0\right) & = &
\frac{-3f \! \left(x_0\right) + 4f \! \left(x_0+h\right) -
f \! \left(x_0+2h\right)}{2h} \\
& + & \frac{h^2}{3}f^{\left(3\right)} \! \left(\xi_0\right) \quad
\textrm{with} \quad \xi_0 \in \left[x_0,x_0+2h\right]. \nonumber
\end{eqnarray}
%
If $x_0$ is the middle point of the sequence,
$\{x_0-h,x_0,x_0+h\}$, then we obtain the familiar ``centered
difference" formula,
\begin{eqnarray}
f^{\prime} \! \left(x_0\right) & = &
\frac{-f \! \left(x_0-h\right) + f \! \left(x_0+h\right)
}{2h} \\
& - & \frac{h^2}{6}f^{\left(3\right)} \! \left(\xi_1\right) \quad
\textrm{with} \quad \xi_1 \in \left[x_0-h,x_0+h\right] \nonumber
\end{eqnarray}
Finally, if $x_0$ is the last of three points,
$\{x_0-2h,x_0-h,x_0\}$, then by replacing $h$ with $-h$ in
Equation (\ref{eqn:derFormL}), we obtain a similar formula for the
derivative at the last point of the sequence.  Note that by the
mean value theorem, with the appropriate choice of $\xi_k$ the
formulas are exact, and in each case in the limit as $h$
approaches zero, are per definition the derivatives at the point
$x_0$. By truncating the remainder terms, we obtain second order
accurate derivatives at each of the three points;
thus for the sequence $\{x_0,x_1,x_2\}$ with even spacing $h$, we
have the following formulas,
\begin{eqnarray}
f^{\prime} \! \left(x_0\right) & \approx &
\frac{-3f \! \left(x_0\right) + 4f \! \left(x_1\right) -
f \! \left(x_2\right)}{2h} \label{eqn:leftForm} \\
f^{\prime} \! \left(x_1\right) & \approx &
\frac{-f \! \left(x_0\right) + f \! \left(x_2\right) }{2h} \label{eqn:centreForm} \\
f^{\prime} \! \left(x_2\right) & \approx &
\frac{f \! \left(x_0\right) - 4f \! \left(x_1\right) + 3f \!
\left(x_2\right)}{2h} \label{eqn:rightForm}
\end{eqnarray}
respectively for the three points.  Figure \ref{fig:threePt} shows
the interpolating polynomial for the three points (a parabola),
whereby the derivatives are denoted by the tangent lines; note
that all three derivatives are of the same interpolating
polynomial.
  This concept is extended to longer sequences of
points, as shown in Figure \ref{fig:fivePt}.  The central formula,
Equation (\ref{eqn:centreForm}), is used everywhere where there
are values to the left and right.  The left and right end-points
use Equations (\ref{eqn:leftForm}) and (\ref{eqn:rightForm}),
respectively.  Clearly, the first two points use the same
interpolating polynomial, and similarly for the last two points.
Obviously, for long sequences of points, keeping track of such
formulas will obscure the structure of the problem at hand.  The
algebra involved in many problems such as surface reconstruction
from gradients is greatly simplified by taking a matrix algebraic
approach to differentiation; specifically, we can write the three
point formulas in matrix form as,
\begin{equation}
\begin{bmatrix} f^{\prime} \! \left(x_0\right) \\ f^{\prime} \! \left(x_1\right) \\ f^{\prime} \! \left(x_2\right)  \end{bmatrix}
\approx \frac{1}{2h}\begin{bmatrix} -3 & 4 & -1 \\ -1 & 0 & 1 \\ 1
& -4 &
3 \end{bmatrix} \begin{bmatrix} f \! \left(x_0\right) \\
f \! \left(x_1\right) \\ f \! \left(x_2\right)
\end{bmatrix}
\end{equation}
Similarly, for the five point sequence shown in Figure
\ref{fig:fivePt}, the appropriate matrix operation to compute the
numerical derivatives is,
\begin{equation}
\begin{bmatrix} f^{\prime} \! \left(x_0\right) \\ f^{\prime} \! \left(x_1\right)
 \\ f^{\prime} \! \left(x_2\right) \\ f^{\prime} \! \left(x_3\right) \\ f^{\prime} \! \left(x_4\right) \end{bmatrix} \approx
\frac{1}{2h}\begin{bmatrix} -3 & 4 & -1 & 0 & 0 \\ -1 & 0 & 1 & 0
& 0 \\ 0 & -1 & 0 & 1 & 0 \\ 0 & 0 & -1 & 0 & 1 \\ 0 & 0 & 1 & -4
& 3
\end{bmatrix}
\begin{bmatrix} f \! \left(x_0\right) \\
f \! \left(x_1\right)
 \\ f \! \left(x_2\right) \\ f \! \left(x_3\right) \\ f \! \left(x_4\right)
 \end{bmatrix}.
\end{equation}
This concept of matrix based numerical differentiation is
fundamental to the methods derived in this paper, since generally,
the numerical differentiation of the discrete function $\V{y} =
f\left(\V{x}\right)$ can be represented and computed by the
matrix-algebraic equation,
\begin{equation}
\V{y}^{\prime} = \M{D}\V{y}.
\end{equation}
Under this premise, we will henceforth omit the $\approx$ under
the contention that
the \textit{numerical derivative} is equal to this relation. \\
  The advantage of this matrix based approach is that without
difficulty, higher order derivative formulas can be used.  For
example, the five point formulas are
\begin{equation}
\begin{bmatrix} f^{\prime} \! \left(x_0\right) \\ f^{\prime} \! \left(x_1\right)
 \\ f^{\prime} \! \left(x_2\right) \\ f^{\prime} \! \left(x_3\right) \\ f^{\prime} \! \left(x_4\right) \end{bmatrix} \approx
\frac{1}{12h}\begin{bmatrix}
-25&48&-36&16&-3\\
-3&-10&18&-6&1\\
1&-8&0&8&-1\\-1&6&-18& 10&3\\3&-16&36&-48&25
\end{bmatrix}
\begin{bmatrix} f \! \left(x_0\right) \\
f \! \left(x_1\right)
 \\ f \! \left(x_2\right) \\ f \! \left(x_3\right) \\ f \! \left(x_4\right)
 \end{bmatrix}
\end{equation}
and are fourth order accurate.  In contrast, methods described in
the literature cannot be extended beyond the use of forward and
backward differences. Besides being only first order accurate,
this unfortunately leads to some inconsistencies.  Specifically,
for a three point sequence, both the formulas,
\begin{equation}
\begin{bmatrix} f^{\prime} \! \left(x_0\right) \\ f^{\prime} \! \left(x_1\right) \\ f^{\prime} \! \left(x_2\right)  \end{bmatrix}
\approx \frac{1}{h}\begin{bmatrix} -1 & 1 & 0 \\ -1 & 1 & 0 \\ 0 &
-1 &
1 \end{bmatrix} \begin{bmatrix} f \! \left(x_0\right) \\
f \! \left(x_1\right) \\ f \! \left(x_2\right)
\end{bmatrix}
\end{equation}
and
\begin{equation}
\begin{bmatrix} f^{\prime} \! \left(x_0\right) \\ f^{\prime} \! \left(x_1\right) \\ f^{\prime} \! \left(x_2\right)  \end{bmatrix}
\approx \frac{1}{h}\begin{bmatrix} -1 & 1 & 0 \\ 0 & -1 & 1 \\ 0 &
-1 &
1 \end{bmatrix} \begin{bmatrix} f \! \left(x_0\right) \\
f \! \left(x_1\right) \\ f \! \left(x_2\right)
\end{bmatrix}
\end{equation}
are first order accurate, however, are obviously different.  When
working with forward and backward differences, at some point in
the sequence one must switch from forward to backward in order
that appropriate formulas are used at the end points.  Commonly,
the sequence forward/central/backward is used resulting in the
operator,
\begin{equation}
\begin{bmatrix} f^{\prime} \! \left(x_0\right) \\ f^{\prime} \! \left(x_1\right) \\ f^{\prime} \! \left(x_2\right)  \end{bmatrix}
\approx \frac{1}{2h}\begin{bmatrix} -2 & 2 & 0 \\ -1 & 0 & 1 \\ 0
& -2 &
2 \end{bmatrix} \begin{bmatrix} f \! \left(x_0\right) \\
f \! \left(x_1\right) \\ f \! \left(x_2\right)
\end{bmatrix}.
\end{equation}
While these are theoretically correct derivative formulas, the are
inconsistent, since the central formula is second order accurate
in contrast to the forward/backward formulas which are only first
order accurate\footnote{This is unfortunately the derivative
approximation used in MATLAB$^{\circledR}$'s \texttt{gradient}
function.}.  Consistency at the endpoints is all the more critical
when considering boundary conditions. Unfortunately, such
discussions of the end points are more often than not completely
avoided in the literature (e.g.,~\cite{horn1986}), and even
altogether incorrect formulas
are used (e.g.,~\cite{agrawal2006}). \\
Throughout this section, for simplicity, it has been assumed that
the data is on evenly spaced points.  However, it is not difficult
to derive the appropriate formulas for arbitrary node spacing,
$h_k$; said formulas have been omitted for clarity.
\subsection{Surface Reconstruction from Gradients}
\label{ssec:GLS}
The novelty of formulating numerical differentiation as a matrix
multiplication is that the partial derivatives of a surface take
the particularly simple form,
\begin{eqnarray}
\pdiff{\M{Z}}{x} & = & \M{Z}\MT{D}_x \\
\pdiff{\M{Z}}{y} & = & \M{D}_y\M{Z}.
\end{eqnarray}
Note that $\M{D}_x$ is defined as a differentiation matrix, as
above, and is transposed to effect differentiation in the
$x$-direction. To address the reconstruction problem, we denote a
measured gradient field, obtained for example via Photometric
Stereo~\cite{woodham1980}, as $\M{\hat{Z}}_x$ and $\M{\hat{Z}}_y$.
The reconstruction problem can then be formulated as finding the
surface $\M{Z}$ such that,
\begin{equation}
\M{\hat{Z}}_x \approx \M{Z}\MT{D}_x \quad \textrm{and} \quad
\M{\hat{Z}}_y \approx \M{D}_y\M{Z}
\end{equation}
where $\approx$ denotes equality in the least-squares sense.
 The global least squares cost function for the
reconstruction of a surface from its gradient
field~\cite{Harker2008c} is therefore written in terms of the
matrix Frobenius norm as,
\begin{equation} \label{eqn:fundCF}
\epsilon(\M{Z}) = \left\| \M{Z}\MT{D}_x - \M{\hat{Z}}_x
\right\|_{\textrm{F}}^2 + \left\| \M{D}_y\M{Z} - \M{\hat{Z}}_y
\right\|_{\textrm{F}}^2,
\end{equation}
which represents the Euclidean distance from the measured gradient
field to the gradient field of an unknown surface $\M{Z}$. From a
mathematical point of view, this cost function can be considered
to be a discrete functional in reference to the calculus of
variations, whereby, it is a function of the unknown function
(surface) $\M{Z}$.  To find the minimum of the cost function, we
differentiate with respect to the matrix\footnote{The relevant
derivative formula can be derived using the trace definition of
the Frobenius norm and the formulas developed by
Sch\"{o}nemann~\cite{schoenemann1985}. } $\M{Z}$, yielding the
effective normal equations of the least-squares problem,
\begin{equation} \label{eqn:normalEq} \MT{D}_y\M{D}_y\M{Z}
+ \M{Z}\MT{D}_x\M{D}_x
- \MT{D}_y\M{\hat{Z}}_y - \M{\hat{Z}}_x\M{D}_x = \M{0}.
\end{equation}
This matrix equation is a set of equations which are linear in the
unknowns $z_{ij}$, and is known as a Sylvester Equation (cf.
Equation (\ref{eqn:sylEq0}) and see also~\cite{Stewart2}). Due to
the fact that differentiation matrices are involved, this equation
is rank-one deficient, which would normally indicate a non-unique
solution~\cite{bartels1972}.  However, it is shown in Section
\ref{sec:SEF} that the solution to this equation is unique up to a
constant of integration, as expected.
\subsection{Numerical Solution of Sylvester Equations}

The most common
approach~\cite{horn1986,karacali2003,horovitz2004,agrawal2006,ng2010,koskulics2012,balzer2012}
to the surface reconstruction from gradients problem proceeds by
``vectorizing" the surface $\M{Z}$. That is, by writing,
\begin{equation}
\V{z} = \begin{bmatrix} \V{z}_1 \\ \vdots \\ \V{z}_n
\end{bmatrix},
\end{equation}
where the $\V{z}_k$ are a column partitioning of $\M{Z}$, i.e.,
\begin{equation}
\M{Z} = \begin{bmatrix} \V{z}_1 & \cdots & \V{z}_n \end{bmatrix}.
\end{equation}
This operation is usually denoted,
\begin{equation}
\V{z} = \mathop{\mathrm{vec}}\left(\M{Z}\right).
\end{equation}
The resulting linear system of equations to be solved is therefore
of the form,
\begin{equation}
\M{A}\V{z}=\V{b},
\end{equation}
where the coefficient matrix $\M{A}$ is $2mn\times mn$. The
relation of these commonly used methods and the approach based on
the Frobenius norm approach proposed in~\cite{Harker2008c} can be
seen
 by means
of applying the ``vec" operator to the cost function in Equation
(\ref{eqn:fundCF}). The cost function in terms of the Frobenius
norm is algebraically equivalent to the standard linear least
squares problem\footnote{Note that this is the vectorization of
Equation (\ref{eqn:fundCF}), and not the typical discretization
found in the literature, where FEM~\cite{lapidus} type
discretizations are used.},
\begin{equation} \label{eqn:normalEqVec}
\begin{bmatrix} \M{D}_x \otimes \M{I}_m \\
\M{I}_n \otimes \M{D}_y \end{bmatrix}
\mathop{\mathrm{vec}}\left(\M{Z}\right) =
\begin{bmatrix} \mathop{\mathrm{vec}}\left(\M{\hat{Z}}_x\right) \\
\mathop{\mathrm{vec}}\left(\M{\hat{Z}}_y\right) \end{bmatrix}
\end{equation}
where $\otimes$ denotes the Kronecker product~\cite{Stewart2}. The
coefficient matrix of this least squares problem is similarly
$2mn\times mn$.  Since the appropriate solution of this problem
requires a Moore-Penrose pseudo-inverse (or its numerical
equivalent), the solution in this manner is necessarily
computationally intensive. The number of floating point operations
to solve this if the coefficient matrix is full is,
\begin{equation}
W_{\textrm{VEC}} = 41m^3n^3,
\end{equation}
and is therefore an $\mathcal{O}\!\left(n^6\right)$ method (cf.
Higham~\cite[Ch.16]{higham}). A much more efficient manner is to
work with the Sylvester Equation directly, i.e., Equation
(\ref{eqn:normalEq}), as proposed
in~\cite{Harker2008c,harker2011}. A common method for solving
Sylvester Equations is that of Bartels and
Stewart~\cite{bartels1972}, and is described in~\cite{Stewart2}. A
generally more efficient solution is the Hessenberg-Schur method
of Golub et al.~\cite{golub1979}. The number of flops (floating
point operations), or work required, to compute the solution using
this method is,
\begin{equation} \label{eqn:workHS}
W_{\textrm{HS}}\left(m,n\right) = \frac{5}{3} m^3 + 10 n^3 + 5
m^2n + \frac{5}{2} mn^2.
\end{equation}
Clearly, the approach via the Sylvester Equation is an
$\mathcal{O}\!\left(n^3\right)$ algorithm, in stark contrast to
the vectorization approach which yields an
$\mathcal{O}\!\left(n^6\right)$ algorithm\footnote{For comparative
purposes, a typical algorithm for the computation of the SVD
requires $W_{\textrm{SVD}}(m,n)=4m^2n+8mn^2+9n^3$ flops.  It could
be argued that the $\mathcal{O}\!\left(n^6\right)$ algorithms
could use sparse methods, but this argument is fruitless.  An
$\mathcal{O}\!\left(n^6\right)$ algorithm has an asymptotic
computation time of $t=\alpha n^6$; sparse methods aim to reduce
the value of $\alpha$, and do not reduce the complexity of the
problem, which is always identical to that of Gaussian
elimination.  To this end, in Section \ref{sec:numTest} we have
computed the solution to one and the same problem using a sparse
$\mathcal{O}\!\left(n^6\right)$ algorithm and the newly proposed
method; the newly proposed method is incomparably faster.  A
further disadvantage of sparse methods is that they typically
terminate before a proper minimizing solution is attained.}; the
significance of this difference can be seen when one considers
that for most real problems $m$ and $n$ will be of the order of
thousands (e.g., megapixel images from Photometric Stereo).  All
of the new regularization methods presented in the following are
shown to fall into the Sylvester Equation framework, and therefore
share in these computational advantages.
%
%
%
\section{Spectral Regularization} \label{sec:specReg}
%
\subsection{Generalized Fourier Series of Discrete Orthogonal Basis Functions}
A Generalized Fourier Series is the series expansion of a function
in terms of a complete set of orthogonal functions,
$\varphi_{k}(x)$, as,
\begin{equation}
f(x) = \sum_{k=0}^{\infty} \alpha_k \varphi_{k}(x).
\end{equation}
Computation of the coefficients, $\alpha_k$, arises from the
weighted least squares approximation of the function by the
series, that is, by minimizing the function,
\begin{equation}
\epsilon(\alpha_0,\ldots,\alpha_{\infty}) = \int_{a}^{b} w(x)
\left(f(x) - \sum_{k=0}^{\infty} \alpha_k \varphi_{k}(x)\right)^2
\textrm{d}x,
\end{equation}
where $w(x)$ is a positive weighting function.  Differentiating
with respect to the $\kth{j}$ coefficient, $\alpha_j$ yields,
\begin{equation}
\pdiff{\epsilon}{\alpha_j} = 2\int_{a}^{b} w(x) \left(f(x) -
\sum_{k=0}^{\infty} \alpha_k \varphi_{k}(x)\right) \varphi_{j}(x)
\, \textrm{d}x,
\end{equation}
whereby equating to zero yields the relation,
\begin{equation}
\int_{a}^{b} w(x) f(x)  \varphi_{j}(x) \textrm{d}x =
\sum_{k=0}^{\infty} \alpha_k \int_{a}^{b} w(x) \varphi_{k}(x)
\varphi_{j}(x) \, \textrm{d}x.
\end{equation}
Thus, if the basis functions satisfy the orthogonality condition,
\begin{equation}
\int_{a}^{b} w(x) \varphi_{k}(x) \varphi_{j}(x) \textrm{d}x =
\kappa_{k}\delta_{kj},
\end{equation}
Then each of the coefficients is given as,
\begin{equation}
\alpha_k = \frac{1}{\kappa_{k}} \int_{a}^{b} w(x) f(x)
\varphi_{k}(x) \, \textrm{d}x.
\end{equation}
When working in the discrete domain, however, the function $f(x)$
is only known at a finite number of points, and hence the
coefficients cannot be computed in this manner.  To this end, we
require so-called discrete orthogonal basis functions, which are
continuous basis functions that are orthogonal over a discrete
measure.  That is to say that the orthogonality condition reads,
\begin{equation}
\int_{a}^{b} w(x) \varphi_{k}(x) \varphi_{j}(x)
\textrm{d}\lambda(x) = \kappa_{k}\delta_{kj},
\end{equation}
where the measure $\lambda(x)$ has the differential\footnote{The
corresponding measure for continuous basis functions is
$\lambda(x) = x$, whereby $\textrm{d}\lambda=\textrm{d}x$.},
\begin{equation}
\textrm{d}\lambda(x) = \sum_{i=1}^{n} \delta(x-x_i)  \,
\textrm{d}x,
\end{equation}
where the $x_i$ are the abscissae where the value of $f(x)$ is
known.  Thus, by the sifting property of the delta function, the
discrete orthogonality condition is,
\begin{equation}
\sum_{i=1}^{n} w(x_i) \varphi_{k}(x_i) \varphi_{j}(x_i) =
\kappa_{k}\delta_{kj}.
\end{equation}
The coefficients for the discrete Generalized Fourier Series are
then given as,
\begin{equation}
\alpha_k = \frac{1}{\kappa_k}\sum_{i=1}^{n} w(x_i) f(x_i)
\varphi(x_i),
\end{equation}
which clearly only requires the values of the function $f(x)$ at
the ordered set of nodes $x_k$, $k=1,\ldots,n$. Discrete basis
functions were originally created as a computationally efficient
solution to the interpolation--approximation
problem\footnote{Failure to make this transition to discrete basis
functions leads to spectral methods akin to finite element
analysis, such as~\cite{balzer2011,balzer2012}, which are
notoriously computationally intensive and impractical.} (e.g., the
Gram polynomials~\cite{gram1883}). To compute the coefficients of
a discrete series, we also only require the values of the basis
functions at the set of nodes, which can be conveniently arranged
and manipulated in matrix form as,
\begin{equation}
\M{B} = \begin{bmatrix}
\varphi_0(x_1) & \varphi_1(x_1) & \cdots & \varphi_{n-1}(x_1) \\
\varphi_0(x_2) & \varphi_1(x_2) & \cdots & \varphi_{n-1}(x_2) \\
\vdots & \vdots & \ddots & \vdots \\
\varphi_0(x_n) & \varphi_1(x_n) & \cdots & \varphi_{n-1}(x_n) \\
\end{bmatrix}
\end{equation}
In this manner, the function $\V{y}=f\!\left(\V{x}\right)$ is then
described as
\begin{equation}
\V{y} = \M{B}\V{\alpha}
\end{equation}
whereby the matrix $\M{B}$ can be viewed as composed column-wise
of vector basis functions,
\begin{equation}
\M{B} = \begin{bmatrix} \V{b}_{0}(\V{x}) & \V{b}_{1}(\V{x}) &
\cdots & \V{b}_{n-1}(\V{x})
\end{bmatrix}
\end{equation}
and satisfies the orthogonality condition,
\begin{equation}
\MT{B}\M{W}\M{B} = \M{I},
\end{equation}
usually with the added condition that the weighting matrix,
$\M{W}$, is positive definite.  Clearly when working with discrete
basis functions, the series expansion of a function is necessarily
finite.  Further, any function (vector) can be represented by such
a series provided that the set of basis functions is complete; the
completeness of the set of basis functions, $\V{b}_k$, entails
that the matrix $\M{B}$ is $n \times n$ and full rank.  To extend
discrete orthogonal basis functions to a 2D domain, we define the
basis functions $\M{B}_x$ and $\M{B}_y$ respectively for the $x$
and $y$ directions.  The surface $\M{Z}$ is then represented as,
\begin{equation}
\M{Z} = \M{B}_y\M{C}\MT{B}_x,
\end{equation}
where the matrix $\M{C}$ represents the generalized Fourier
coefficients.  If the basis functions satisfy the orthogonality
conditions,
\begin{eqnarray}
\MT{B}_x\M{W}_x\M{B}_x & = & \M{I}_n \\
\MT{B}_y\M{W}_y\M{B}_y & = & \M{I}_m
\end{eqnarray}
then the minimizing coefficients of the weighted least squares
cost function (weighted Frobenius norm),
\begin{equation}
\epsilon(\M{C}) = \left\| \M{W}_y^{\frac{1}{2}}\left( \M{Z} -
\M{B}_y\M{C}\MT{B}_x \right) \M{W}_x^{\frac{1}{2}}
\right\|_{\textrm{F}}^2
\end{equation}
are the generalized Fourier coefficients,
\begin{equation} \label{eqn:gfsCoeffs}
\M{C} = \MT{B}_y\M{W}_y\M{Z}\M{W}_x\M{B}_x
\end{equation}
Some examples of functions which can be considered to be
generalized Fourier series in the discrete sense are Gram
Polynomials~\cite{oleary2008b}, Cosine/Sine Functions (e.g. DCT),
Fourier Basis Functions, Haar Functions~\cite{haar1910}, Hartley
Transform Basis Functions~\cite{bracewell}, etc.
\subsection{Surface Reconstruction with GFS}
According to Equation (\ref{eqn:gfsCoeffs}) we can represent a
reconstructed surface exactly with any complete set of basis
functions.  Hence using the complete set of basis functions should
have no influence on the surface reconstruction other than to
increase the computational load.  On the other hand, if we use an
incomplete (or truncated) set of basis functions, then we can
effectively incorporate band-pass filtering into the least squares
solution. Specifically, we represent $\M{Z}$ as the truncated
generalized Fourier series,
\begin{equation} \label{eqn:gfsSurf}
\M{Z} = \M{B}_y\M{C}\MT{B}_x,
\end{equation}
where $\M{B}_y$ is a set of $p$ basis functions on $m$ nodes with
$p < m$, matrix $\M{B}_x$ is a set of $q$ basis functions on $n$
nodes with $q < n$, and the matrix $\M{C}$ is a $p \times q$
matrix of generalized Fourier coefficients.  For simplicity, we
assume the functions are orthogonal with respect to an identity
weighting, that is, the basis functions are orthonormal (The
weighted least squares solution is postponed until Section
\ref{sec:WLS}). Thus the least squares cost function is obtained
by substituting the surface representation in Equation
(\ref{eqn:gfsSurf}) into Equation (\ref{eqn:fundCF}), i.e.,
\begin{equation}
\epsilon(\M{C}) = \left\| \M{B}_y\M{C}\MT{B}_x\MT{D}_x -
\M{\hat{Z}}_x \right\|_{\textrm{F}}^2 + \left\|
\M{D}_y\M{B}_y\M{C}\MT{B}_x - \M{\hat{Z}}_y
\right\|_{\textrm{F}}^2
\end{equation}
Differentiating with respect to the unknown coefficients, $\M{C}$,
yields the effective normal equations,
\begin{eqnarray} \label{eqn:normalEqSpec}
\MT{B}_y\MT{D}_y\M{D}_y\M{B}_y\M{C} & + &
\M{C}\MT{B}_x\MT{D}_x\M{D}_x\M{B}_x \nonumber \\
& - & \MT{B}_y\left( \MT{D}_y \M{\hat{Z}}_y + \M{\hat{Z}}_x\M{D}_x
\right)\M{B}_x = \M{0},
\end{eqnarray}
which is a $p \times q$ Sylvester Equation in the unknown
coefficients $\M{C}$.
%
%
Since a surface defined as a function of its spectral coefficients
has, per definition, an integrable gradient field, its gradient
field spans a subspace of all integrable gradient fields.  That
is, its gradient spans a band limited subspace of the integrable
gradient fields.
\subsection{Computational Aspects}
The solution of the surface reconstruction problem is obtained by
solving the Sylvester Equation (\ref{eqn:normalEqSpec}) and
back-substitution of the coefficients into Equation
(\ref{eqn:gfsSurf}). Using the method of Golub et
al.~\cite{golub1979}, the computational work required is given by
Equation (\ref{eqn:workHS});
%
%
however, in the case of Spectral Regularization the solution of
the Sylvester Equation is more efficient since an $m\times n$
equation is reduced to a $p\times q$ equation. That is, if $p$ and
$q$ are some fraction of $m$ and $n$ of the form
\begin{equation}
p = \frac{m}{2^k} \qquad \textrm{and} \qquad q = \frac{n}{2^k},
\end{equation}
then the number of flops to solve the corresponding Sylvester
Equation is,
\begin{equation}
W\left(p,q\right) = \frac{1}{2^{3k}}W\left(m,n\right)
\end{equation}
Hence, by using only half of the basis functions $(k=1)$, the work
is reduced to $\frac{1}{8}$ of the full problem.  By using one
quarter of the basis functions $(k=2)$, the work is reduced to
$\frac{1}{64}$, etc.
%
%
\section{Tikhonov Regularization} \label{sec:tikhReg}

Tikhonov regularization\footnote{Tikhonov Regularization goes
under a number of other pseudonyms: Tikhonov-Phillips
regularization, ridge regression, damped least squares,
etc.}~\cite{Engl} over a 1D domain amounts to finding the function
$\V{y}$ which minimizes the functional,
\begin{equation}
\epsilon\left(\V{y}\right) = \left\| \M{A}\V{y} - \V{b}
\right\|_2^2 + \lambda^2 \left\| \M{L}\left(\V{y}-\V{y}_0\right)
\right\|_2^2.
\end{equation}
Whereas the former term is a typical least-squares cost function,
the latter acts as a ``penalty" term. The function $\V{y}_0$ is an
\textit{a priori} estimate of the unknown function; if nothing is
known about the function, then $\V{y}_0$ is the zero vector. Hence
the penalty term is a (weighted) measure of the deviation from the
\textit{a priori} estimate.  Thus if $\M{L} = \M{I}$, then the
penalty term is the Euclidean deviation of the solution from the
\textit{a priori} estimate.  The constant $\lambda$ is the
regularization parameter, which is positive and is assumed to be
fixed for the optimization process; it essentially shifts the
priority between the least squares residual and the regularization
term. To find the minimum of the functional, it is differentiated
with respect to $\V{y}$, where upon rearranging yields,
\begin{equation} \label{eqn:tikhNE}
\left(\MT{A}\M{A} + \lambda^2\MT{L}\M{L}\right)\V{y} = \MT{A}\V{b}
+ \lambda^2\MT{L}\M{L}\V{y}_0,
\end{equation}
which are the corresponding normal equations.
 Equation (\ref{eqn:tikhNE}) is a
\textit{necessary}, but not \textit{sufficient}, condition that
$\V{y}$ is a constrained minimizer of the least-squares cost
function. For this reason, regularization is typically associated
with Lagrange Multipliers; however, it is actually much more
closely related to the Levenberg-Marquardt
algorithm~\cite{Marquardt1963}; specifically, each step of the
Levenberg-Marquardt algorithm is a Tikhonov-reg\-ulari\-zed
Gauss-Newton step. For a fixed $\lambda$ the minimizing solution
can be computed directly via the Moore-Penrose pseudo-inverse.
The Tikhonov regularization problem is said to be in standard form
if the smoothing operator $\M{L}=\M{I}$ and $\V{y}_0=\V{0}$.
\subsection{Tikhonov Regularization for the Reconstruction Problem}
For the problem of surface reconstruction from gradients, the
functional to be minimized depends on a 2D surface, as opposed to
a vector; hence, the common approach to Tikhonov regularization
does not apply directly.  To derive the appropriate functional for
a 2D domain, we begin with the least squares cost function of
Equation (\ref{eqn:fundCF}).
We define the matrices $\M{L}_x$ and $\M{L}_y$ as general
``smoothing" operators in the $x$- and $y$-directions,
respectively; thus, the functional for Tikhonov regularization in
its most general form reads,
\begin{eqnarray}
\epsilon\left(\M{Z}\right) & = & \left\| \M{D}_y\M{Z} -
\M{\hat{Z}}_y \right\|_{\textrm{F}}^2 + \left\| \M{Z}\MT{D}_x -
\M{\hat{Z}}_x \right\|_{\textrm{F}}^2 \nonumber \\
& + &
\mu^2 \left\| \M{L}_y \left( \M{Z} - \M{Z}_0 \right)
\right\|_{\textrm{F}}^2 +
\lambda^2 \left\| \left( \M{Z} - \M{Z}_0 \right)
\M{L}_x^{\textrm{T}} \right\|_{\textrm{F}}^2.
\end{eqnarray}
where $\M{Z}_0$ is an \textit{a-priori} estimate of the surface.
The estimate $\M{Z}_0$ is not necessary, since
 we may assume the surface is nearly flat, i.e., $\M{Z}_0 = \M{0}$.
 However, the fact that we can incorporate this \textit{a priori}
 estimate into the algorithm may have substantial consequences for
 applied Photometric Stereo.
For the sake of generality, we have introduced a second
regularization parameter, $\mu$, which may be of use if the $x$-
and $y$-derivatives have different noise properties\footnote{For a
treatment of multiple regularization parameters in 1D problems,
see~\cite{belge2002}.}.
Differentiating the cost function with respect to $\M{Z}$ yields
the corresponding normal equations of this functional, i.e.,
\begin{eqnarray} \label{eqn:tikhSE}
\left( \MT{D}_y\M{D}_y + \mu^2\MT{L}_y\M{L}_y\right)\M{Z}
& + & \M{Z}\left( \MT{D}_x\M{D}_x +
\lambda^2\MT{L}_x\M{L}_x\right)
\nonumber \\
& - & \MT{D}_y\M{\hat{Z}}_y - \M{\hat{Z}}_x\M{D}_x \nonumber \\
& - & \mu^2\MT{L}_y\M{L}_y\M{Z}_0 \nonumber \\
& - & \M{Z}_0\left(\lambda^2\MT{L}_x\M{L}_x\right) = \M{0}.
\end{eqnarray}
This is again a Sylvester Equation in the unknown surface $\M{Z}$.
For the surface reconstruction problem, the standard form of the
Tikhonov regularization problem corresponds to $\M{L}_x =
\M{I}_n$, $\M{L}_y = \M{I}_m$, and $\M{Z}_0=\M{0}$, in which case
it suffices to consider only $\mu = \lambda$.
\subsection{Regularization Terms}
In the context of a 2D reconstruction problem, the following
regularization terms were derived in their matrix form
in~\cite{harker2011}.  For completeness, the results are merely
summarized here. The most basic regularization term is a bound on
the norm of the solution, in this case,
\begin{equation}
\rho\left(\M{Z}\right) = \left\| \M{Z} \right\|_{\textrm{F}}^2,
\end{equation}
which is effectively a degree-0 regularization term, and
corresponds to the Tikhonov problem in its standard form. It is
written equivalently as,
\begin{equation}
\rho\left(\M{Z}\right) = \frac{1}{2} \left( \left\| \M{I}_m\M{Z}
\right\|_{\textrm{F}}^2 + \left\| \M{Z}\M{I}_n
\right\|_{\textrm{F}}^2 \right),
\end{equation}
such that it corresponds to the Sylvester Equation
(\ref{eqn:tikhSE}).
A degree-1 regularization term bounds the magnitude of the maximum
directional derivative, or the overall steepness of the
reconstructed surface, i.e.,
\begin{equation}
\rho\left(\M{Z}\right) =
\left\| \M{D}_y\M{Z} \right\|_{\textrm{F}}^2 + \left\|
\M{Z}\MT{D}_x \right\|_{\textrm{F}}^2.
\end{equation}
Finally, the regularization term,
\begin{equation}
\rho\left(\M{Z}\right) = \left\| \M{D}_y^2\M{Z}
\right\|_{\textrm{F}}^2 +
 \left\| \M{Z}\left(\M{D}_x^2\right)^{\textrm{T}}
 \right\|_{\textrm{F}}^2,
\end{equation}
is a degree-2 regularization term, and bounds the mean curvature
of the surface.
\subsection{Influence of the Regularization Parameter}
To derive an effective algorithm for determining the
regularization parameter, as well as to characterize the effect of
the regularization parameter on the solution, we look at the
Tikhonov problem in its standard form.  Firstly, we denote the
SVDs of the $x$- and $y$-derivative operators as
\begin{eqnarray}
\M{D}_x & = & \M{U}_x\M{S}_x\MT{V}_x = \sum_{i=1}^{n}\alpha_i\V{u}_i\VT{v}_i \\
\M{D}_y & = & \M{U}_y\M{S}_y\MT{V}_y =
\sum_{j=1}^{m}\beta_j\V{m}_j\VT{w}_j
\end{eqnarray}
By substitution of these relations, the normal equations for
Tikhonov regularization in standard form can be written as,
\begin{eqnarray}
\left(\M{V}_y\M{S}_y^2\MT{V}_y + \lambda^2\M{I}_m\right)\M{Z} & +
& \M{Z}\left(\M{V}_x\M{S}_x^2\MT{V}_x + \lambda^2\M{I}_n\right)
\nonumber \\ - \M{V}_y\M{S}_y\MT{U}_y \M{\hat{Z}}_y & - &
\M{\hat{Z}}_x\M{U}_x\M{S}_x\MT{V}_x = \M{0}.
\end{eqnarray}
Pre-multiplying by $\MT{V}_y$ and post-multiplying by $\M{V}_x$
yields,
\begin{eqnarray}
\left(\M{S}_y^2\MT{V}_y + \lambda^2\MT{V}_y\right)\M{Z}\M{V}_x & +
&
\MT{V}_y\M{Z}\left(\M{V}_x\M{S}_x^2 + \lambda^2\M{V}_x\right) \nonumber \\
- \M{S}_y\MT{U}_y \M{\hat{Z}}_y\M{V}_x & - &
\MT{V}_y\M{\hat{Z}}_x\M{U}_x\M{S}_x = \M{0}.
\end{eqnarray}
If we make the following substitution,
\begin{equation} \label{eqn:ZasFofV}
\M{M} = \MT{V}_y\M{Z}\M{V}_x,
\end{equation}
then the matrix $\M{M}$ represents the generalized Fourier Series
coefficients of the surface $\M{Z}$ with respect to the singular
vectors $\M{V}_x$ and $\M{V}_y$.  Making the further
substitutions,
\begin{eqnarray}
\M{P} & = & \MT{U}_y \M{\hat{Z}}_y\M{V}_x \\
\M{Q} & = & \MT{V}_y \M{\hat{Z}}_x\M{U}_x
\end{eqnarray}
then the normal equations read,
\begin{equation}
\left(\M{S}_y^2 + \lambda^2\M{I}_m\right)\M{M} +
\M{M}\left(\M{S}_x^2 + \lambda^2\M{I}_n\right) - \M{S}_y\M{P} -
\M{Q}\M{S}_x = \M{0}.
\end{equation}
Since all the pertinent coefficient matrices are diagonal, this
equation, taken element-wise, reads,
\begin{equation}
\left(\beta_{i}^2 + \lambda^2\right)m_{ij} +
m_{ij}\left(\alpha_{j}^2 + \lambda^2\right) - \beta_{i}p_{ij} -
\alpha_{j}q_{ij} = 0,
\end{equation}
for each $i=1,\ldots,m$, $j=1,\ldots,n$.  Thus, the entries of
$\M{M}$ can be solved for as a function of $\lambda$ as,
\begin{equation}
m_{ij}\left(\lambda\right) = \frac{\beta_i p_{ij} + \alpha_j
q_{ij} }{ \alpha_j^2 + \beta_i^2 + 2\lambda^2}
\end{equation}
which determines the coefficients in terms of $\lambda$.  To
demonstrate the influence of the regularization parameter,
$\lambda$, this relation can be rewritten in the form,
\begin{equation} \label{eqn:solWithFF}
m_{ij}\left(\lambda\right) =
f_{ij}\left(\lambda\right)\frac{\beta_i p_{ij} + \alpha_j q_{ij}
}{ \alpha_j^2 + \beta_i^2 }
\end{equation}
where the terms
\begin{equation}
f_{ij}\left(\lambda\right) = \frac{ \alpha_j^2 + \beta_i^2 }{
\alpha_j^2 + \beta_i^2 + 2\lambda^2}
\end{equation}
can be considered to be filter factors~\cite{hansen1993}, which
range from $1$ to $0$ (corresponding respectively to $\lambda=0$
and $\lambda\rightarrow\infty$). Clearly, when the filter factors
are all one (i.e., $f_{ij}(0)=1$), Equation (\ref{eqn:solWithFF})
represents the least-squares solution (GLS) to the reconstruction
problem.  The values,
\begin{equation}
\mu_{ij}^2 = \alpha_j^2 + \beta_i^2,
\end{equation}
are the eigenvalues of the Sylvester Operator~\cite{Stewart2}, and
thus writing the filter factors as,
\begin{equation}
f_{ij}\left(\lambda\right) = \frac{ \mu_{ij}^2 }{ \mu_{ij}^2 +
2\lambda^2},
\end{equation}
shows that Tikhonov regularization for the reconstruction problem
in Sylvester Equation form, has essentially the same structure as
standard 1-dim\-ension\-al domain Tikhonov regularization
problems, cf.~\cite{hansen1993}.
 As can be seen from the
filter factors, the regularization parameter $\lambda$ inversely
weights the coefficients $m_{ij}$. This influence of $\lambda$ on
the reconstructed surface can be seen from Equation
(\ref{eqn:ZasFofV}), since the reconstructed surface can be
written as,
\begin{equation}
\M{Z}\left(\lambda\right) = \sum_{i=1}^{m}\sum_{j=1}^{n}
m_{ij}\left(\lambda\right) \V{w}_{i} \VT{v}_{j},
\end{equation}
The reconstructed surface is therefore sum of rank-1 matrices each
weighted by the coefficients $m_{ij}$, which are functions of the
regularization parameter $\lambda$.  Clearly, the parameter
$\lambda$ has a larger influence on basis functions corresponding
to small singular values; specifically, basis functions associated
with small singular values are largely suppressed.  Note that this
analysis constitutes an algorithm for solving symmetric
rank-deficient Sylvester Equations; however, here it is used as a
means of effectively determining the regularization parameter
$\lambda$.
\subsection{Selection of the Regularization Parameter}
The L-curve is a plot of $\left( \rho\left(\lambda\right),
\eta\left(\lambda\right) \right)$ where,
$\rho^2\left(\lambda\right)$ is the least-squares cost function,
\begin{equation}
\rho^2\left(\lambda\right) = \left\|
\M{Z}\left(\lambda\right)\MT{D}_x - \M{\hat{Z}}_x
\right\|_{\textrm{F}}^2 + \left\| \M{D}_y\M{Z}\left(\lambda\right)
- \M{\hat{Z}}_y \right\|_{\textrm{F}}^2
\end{equation}
and $\eta^2\left(\lambda\right)$ is the regularization term (in
standard form),
\begin{equation}
\eta^2\left(\lambda\right) = \left\| \M{Z}\left(\lambda\right)
\right\|_{\textrm{F}}^2
\end{equation}
Thus, the L-curve is a visualization of the interplay of the least
squares residual, and the magnitude of the regularization term.
Once the singular value decompositions of the derivative matrices
are computed, points on the L-curve can be computed as,
\begin{equation}
\rho^2\left(\lambda\right) = \left\|
\M{M}\left(\lambda\right)\M{S}_x - \M{Q} \right\|_{\textrm{F}}^2 +
\left\| \M{S}_y\M{M}\left(\lambda\right) - \M{P}
\right\|_{\textrm{F}}^2
\end{equation}
and
\begin{equation}
\eta^2\left(\lambda\right) = \left\| \M{M}\left(\lambda\right)
\right\|_{\textrm{F}}^2
\end{equation}
due to the invariance of the Frobenius norm under orthonormal
transformation.  The computational cost of these evaluations are
relatively small in comparison to the computation of the singular
value decomposition, due to the fact that $\M{S}_x$ and $\M{S}_y$
are diagonal. The result is that several points on the L-curve can
be computed to determine an appropriate regularization parameter.
Further analysis of regularization parameter selection is beyond
the scope of this paper; the reader is referred
to~\cite{Engl,GolubMeurant}.
%
%
\section{Constraining Solutions by Known Boundary Conditions}
\label{sec:BCs}
%
Boundary conditions are usually imposed as auxiliary conditions to
partial differential equations to ensure the existence of a unique
solution.  For the surface reconstruction from gradients problem,
they can be imposed to constrain the solution, which effects a
form of regularization on the reconstructed surface.  Dirichlet
Boundary Conditions specify the value of the integral surface on
the domain boundary, and can have highly effective regularizing
effects (See Section \ref{sec:numTest}).  Neumann Boundary
Conditions specify the value of the normal derivative of the
integral surface on the domain boundary, and hence has a similar
effect to first order Tikhonov Regularization.  In the following,
we derive the solution to the reconstruction problem with
arbitrary Dirichlet Conditions.
\subsection{Dirichlet Boundary Conditions}
Dirichlet boundary conditions specify the value of the function
(in this case the height of the surface) on the boundary of the
domain.  Using a matrix based approach, we start by parameterizing
a surface with fixed height on the boundary. This can be
accomplished with permutation matrices\footnote{It should be noted
that permutation matrices should not be implemented explicitly, as
they serve the same function as re-indexing the rows and columns
of matrices~\cite[pp.109-110]{Golub}.}; specifically,
\begin{equation}
\M{Z} = \M{P}_m\M{Z}_\textrm{I}\MT{P}_n + \M{Z}_\textrm{B},
\end{equation}
where $\M{P}_m$ and $\M{P}_n$ are the (orthonormal) permutation
matrices,
\begin{equation}
\M{P}_m = \begin{bmatrix} \VT{0} \\ \M{I}_{m-2} \\ \VT{0}
\end{bmatrix}
\qquad \textrm{and} \qquad
\M{P}_n = \begin{bmatrix} \VT{0} \\ \M{I}_{n-2} \\ \VT{0}
\end{bmatrix}.
\end{equation}
The matrix $\M{Z}_\textrm{I}$ is the
$\left(m-2\right)\times\left(n-2\right)$ matrix of the unknown
interior values of the surface and $\M{Z}_\textrm{B}$ is the
$m\times n$ matrix specifying the boundary values\footnote{Here it
is implied that the interior values of $\M{Z}_\textrm{B}$ are
zero, however, this is not necessary.  If some interior values of
$\M{Z}_\textrm{B}$ are non-zero, then $\M{Z}_\textrm{I}$ simply
represents the deviation from this surface.  For example, if
$\M{Z}_\textrm{B}$ specifies a parabolic surface, then
$\M{Z}_\textrm{I}$ would represent the deviation of the internal
portion from this parabolic surface.}.  Substituting this
parametrization into the cost function for the surface
reconstruction problem, i.e., Equation (\ref{eqn:fundCF}), yields,
\begin{eqnarray}
\epsilon(\M{Z}_\textrm{I}) & = & \left\|
\left(\M{P}_m\M{Z}_\textrm{I}\MT{P}_n +
\M{Z}_\textrm{B}\right)\MT{D}_x - \M{\hat{Z}}_x
\right\|_{\textrm{F}}^2 \nonumber \\
& + & \left\| \M{D}_y\left(\M{P}_m\M{Z}_\textrm{I}\MT{P}_n +
\M{Z}_\textrm{B}\right) - \M{\hat{Z}}_y \right\|_{\textrm{F}}^2
\end{eqnarray}
Differentiating the functional with respect to $\M{Z}_\textrm{I}$
yields the effective normal equations,
\begin{eqnarray}
\MT{P}_m\MT{D}_y\M{D}_y\M{P}_m \M{Z}_\textrm{I} & + &
\M{Z}_\textrm{I} \MT{P}_n\MT{D}_x\M{D}_x\M{P}_n \nonumber \\
& - & \MT{P}_m\left( \M{\hat{Z}}_x - \M{Z}_\textrm{B}\MT{D}_x
\right)\M{D}_x\M{P}_n \nonumber \\
& - & \MT{P}_m\MT{D}_y\left( \M{\hat{Z}}_y -
\M{D}_y\M{Z}_\textrm{B} \right)\M{P}_n = \M{0},
\label{eqn:dirichletNE}
\end{eqnarray}
which is an $\left(m-2\right)\times\left(n-2\right)$ Sylvester
Equation in the unknown interior portion of the surface
$\M{Z}_\textrm{I}$.
%
\section{Weighted Least Squares Solutions} \label{sec:WLS}
%
Weighted least squares is an important extension of the standard
least squares problem when measurement errors behave according to
heteroscedastic Gaussian distributions. The maximum likelihood
cost function is the standard least squares cost function modified
by the inverse square root of the covariance matrix of the errors
in measuring $\V{b}$, i.e.,
\begin{equation}
\V{y}_\textrm{M} = \min_{\V{y}}\left\|
\M{\Lambda}^{-\frac{1}{2}}\left( \M{A}\V{y} - \V{b} \right)
\right\|_2^2
\end{equation}
This cost function corresponds to the Mahalanobis distance between
$\M{A}\V{y}$ and $\V{b}$, whereby, minimization proceeds by
differentiating with respect to $\V{y}$ and solving the
corresponding normal equations.
\subsection{Mahalanobis Distance between two Gradient Fields}
If we assume that the errors in measuring a gradient field are
non-uniform and covariant, then the measured gradient field is
related to the true gradient field by,
\begin{eqnarray} \label{eqn:mahalEqns}
\M{\hat{Z}}_x & = & \M{Z}_x +
\M{\Lambda}_{xy}^{\frac{1}{2}}\M{\Delta}_{x}\M{\Lambda}_{xx}^{\frac{1}{2}} \\
\M{\hat{Z}}_y & = & \M{Z}_y +
\M{\Lambda}_{yy}^{\frac{1}{2}}\M{\Delta}_{y}\M{\Lambda}_{yx}^{\frac{1}{2}}
\end{eqnarray}
where $\M{\Delta}_x$ and $\M{\Delta}_y$ are matrices of i.i.d.
Gaussian random variables; the covariance matrices,
$\M{\Lambda}_{uv}$, denote the covariance of the $u$-derivative in
the $v$-direction, and for simplicity the matrix square root is
specifically its symmetric square root. In this case, the measured
gradient is related to the true gradient in terms of the
Mahalanobis distance. The minimum Mahalanobis distance is then
characterized by the least squares minimization of the term,
\begin{equation}
\left\| \M{\Delta}_x \right\|_{\textrm{F}}^2 + \left\|
\M{\Delta}_y \right\|_{\textrm{F}}^2
\end{equation}
Therefore, the Mahalanobis distance between two gradient fields is
given as,
\begin{eqnarray}
\epsilon_{\textrm{M}} & = & \left\|
\M{\Lambda}_{xy}^{-\frac{1}{2}} \left( \M{Z}_x - \M{\hat{Z}}_x
\right) \M{\Lambda}_{xx}^{-\frac{1}{2}} \right\|_{\textrm{F}}^2
\nonumber \\ & + &
 \left\| \M{\Lambda}_{yy}^{-\frac{1}{2}} \left( \M{Z}_y -
\M{\hat{Z}}_y
 \right)\M{\Lambda}_{yx}^{-\frac{1}{2}}
\right\|_{\textrm{F}}^2
\end{eqnarray}
which can be considered to be a weighted Frobenius norm, as
proposed in~\cite{higham2002}.
\subsection{Weighted Solution for Surface Reconstruction}
Given the expression for the Mahalanobis between two gradient
fields, the weighted least squares surface reconstruction from
gradients can be posed as the minimization of the functional,
\begin{eqnarray}
\epsilon\left(\M{Z}\right) & = & \left\|
\M{\Lambda}_{xy}^{-\frac{1}{2}} \left( \M{Z}\MT{D}_x -
\M{\hat{Z}}_x \right) \M{\Lambda}_{xx}^{-\frac{1}{2}}
\right\|_{\textrm{F}}^2 \nonumber \\
 & + & \left\| \M{\Lambda}_{yy}^{-\frac{1}{2}} \left( \M{D}_y\M{Z} - \M{\hat{Z}}_y
 \right)\M{\Lambda}_{yx}^{-\frac{1}{2}}
\right\|_{\textrm{F}}^2
\end{eqnarray}
Differentiating the functional with respect to $\M{Z}$ yields the
effective normal equations, i.e.,
\begin{eqnarray}
\MT{D}_y\M{\Lambda}_{yy}^{-1}\M{D}_y \M{Z}\M{\Lambda}_{yx}^{-1}
& + &
\M{\Lambda}_{xy}^{-1}\M{Z}\MT{D}_x\M{\Lambda}_{xx}^{-1}\M{D}_x \\
& - &
\MT{D}_y\M{\Lambda}_{yy}^{-1}\M{\hat{Z}}_y\M{\Lambda}_{yx}^{-1}
- \M{\Lambda}_{xy}^{-1}\M{\hat{Z}}_x\M{\Lambda}_{xx}^{-1}\M{D}_x =
\M{0}, \nonumber
\end{eqnarray}
which is not immediately a Sylvester Equation.  However, by
pre-multiplying by $\M{\Lambda}_{xy}^{\frac{1}{2}}$,
post-multiplying by $\M{\Lambda}_{yx}^{\frac{1}{2}}$, and
inserting the expressions,
\begin{equation}
\M{\Lambda}_{xy}^{\frac{1}{2}}\M{\Lambda}_{xy}^{-\frac{1}{2}} =
\M{I}
\qquad \textrm{and} \qquad
\M{\Lambda}_{yx}^{-\frac{1}{2}}\M{\Lambda}_{yx}^{\frac{1}{2}} =
\M{I},
\end{equation}
we yield the symmetric form\footnote{Under the assumption that the
$\M{\Lambda}_{uv}$ are full-rank, these transformations do not
alter the solution to the Sylvester Equation.},
\begin{eqnarray} \label{eqn:weightedSE}
\left( \M{\Lambda}_{xy}^{\frac{1}{2}}
\MT{D}_y\M{\Lambda}_{yy}^{-1}\M{D}_y
\M{\Lambda}_{xy}^{\frac{1}{2}} \right) \M{Z}_\textrm{w}
& + & \M{Z}_\textrm{w} \left( \M{\Lambda}_{yx}^{\frac{1}{2}}
\MT{D}_x\M{\Lambda}_{xx}^{-1}\M{D}_x
\M{\Lambda}_{yx}^{\frac{1}{2}} \right) \nonumber \\
& - & \M{\Lambda}_{xy}^{\frac{1}{2}} \MT{D}_y\M{\Lambda}_{yy}^{-1}
\M{\hat{Z}}_y \M{\Lambda}_{yx}^{-\frac{1}{2}} \\
& - &
\M{\Lambda}_{xy}^{-\frac{1}{2}}\M{\hat{Z}}_x\M{\Lambda}_{xx}^{-1}\M{D}_x
\M{\Lambda}_{xy}^{\frac{1}{2}} = \M{0} \nonumber
\end{eqnarray}
which is a Sylvester Equation in the unknown ``weighted" surface
$\M{Z}_\textrm{w}$, where,
\begin{equation}
\M{Z}_\textrm{w} =
\M{\Lambda}_{xy}^{-\frac{1}{2}}\M{Z}\M{\Lambda}_{yx}^{-\frac{1}{2}}.
\end{equation}
The weighted Sylvester Equation in Equation (\ref{eqn:weightedSE})
can be clarified with some simplified notation, i.e., it can be
written in the form,
\begin{equation} \label{eqn:weightedSEmod}
\MT{A}\M{A}\M{Z}_\textrm{w} + \M{Z}_\textrm{w}\MT{B}\M{B} -
\MT{A}\M{F} - \M{G}\M{B} = \M{0}
\end{equation}
where $\M{A}$ and $\M{B}$ can be considered to be weighted
differential operators,
\begin{eqnarray}
\M{A} & = & \M{\Lambda}_{yy}^{-\frac{1}{2}}\M{D}_y
\M{\Lambda}_{xy}^{\frac{1}{2}} \\
\M{B} & = & \M{\Lambda}_{xx}^{-\frac{1}{2}}\M{D}_x
\M{\Lambda}_{yx}^{\frac{1}{2}}
\end{eqnarray}
and $\M{F}$ and $\M{G}$ can be considered to be the weighted
gradient field,
\begin{eqnarray}
\M{F} & = & \M{\Lambda}_{yy}^{-\frac{1}{2}}\M{\hat{Z}}_y
\M{\Lambda}_{yx}^{-\frac{1}{2}} \\
\M{G} & = & \M{\Lambda}_{xy}^{-\frac{1}{2}}\M{\hat{Z}}_x
\M{\Lambda}_{xx}^{-\frac{1}{2}}
\end{eqnarray}
Upon solving Equation (\ref{eqn:weightedSE}), or equivalently,
Equation (\ref{eqn:weightedSEmod}), for $\M{Z}_\textrm{w}$ the
weighted least squares solution is obtained as,
\begin{equation}
\M{Z} =
\M{\Lambda}_{xy}^{\frac{1}{2}}\M{Z}_\textrm{w}\M{\Lambda}_{yx}^{\frac{1}{2}}.
\end{equation}
%
\section{Computational Framework} \label{sec:SEF}
%
All of the methods for surface reconstruction from gradients
presented in this paper have been shown to be solved by means of a
particular Sylvester Equation.  In a more unifying sense, all of
the systems of normal equations have the same form, and hence, we
can place all methods within a common computational framework.
Specifically, all methods presented here have normal equations of
the form,
\begin{equation} \label{eqn:genNE}
\MT{A}\M{A}\M{\Phi} + \M{\Phi}\MT{B}\M{B} - \MT{A}\M{F} -
\M{G}\M{B} = \M{0},
\end{equation}
whereby $\M{\Phi}$ represents the unknown parameters of $\M{Z}$
such that $\M{Z} = f\left(\M{\Phi}\right)$, the matrices $\M{A}$,
$\M{B}$ depend on model parameters (i.e., the differential
operators), and $\M{F}$ and $\M{G}$ depend on both model
parameters as well as the measured data. For example, the general
Tikhonov normal equations, Equation (\ref{eqn:tikhNE}), can be
written in this manner with, $\M{\Phi} = \M{Z}$, and
\begin{equation}
\M{A} = \begin{bmatrix} \M{D}_y \\ \mu \M{L}_y \end{bmatrix},
\qquad \M{B} = \begin{bmatrix} \M{D}_x \\ \lambda \M{L}_x
\end{bmatrix}
\end{equation}
and
\begin{equation}
\M{F} = \begin{bmatrix} \M{\hat{Z}}_y \\ \mu \M{L}_y \M{Z}_0
\end{bmatrix},
\qquad \M{G} = \begin{bmatrix} \M{\hat{Z}}_x & \lambda \M{Z}_0
\MT{L}_x
\end{bmatrix}.
\end{equation}
Table \ref{tab:sylEqnKey} contains a summary of all methods
presented in this paper, and the appropriate coefficient matrices
such that they fit into the Sylvester Equation framework.
\subsection{Solution of Symmetric Semi-Definite Sylvester Equations}
Most of the Sylvester Equations presented here are rank deficient,
and hence special care must be taken in solving them. The null
spaces of all the Sylvester Operators are known \textit{a priori},
and hence valuable computation time need not be wasted by
computing them (e.g., via the SVD).  Specifically, for a simple
differential operator, the derivative of a constant function must
vanish, and hence the null space is fully described as,
\begin{equation}
\M{D}\V{1} = \V{0}.
\end{equation}
Clearly, for an appropriately defined differential operator, this
relation must hold, in addition to the fact that the null space be
of dimension one\footnote{This is a good test that a proper
differential operator has been proposed, since there are several
examples of operators in the literature which do not satisfy these
properties, and are hence not differential operators \textit{per
se} (e.g., \cite{agrawal2006}).}.  According to the proposed
framework, some of the Sylvester Equations use ``modified"
differential operators, $\M{A}$ and $\M{B}$ which are of the
general form,
\begin{equation}
\M{A} = \M{M}\M{D}_y\M{N},
\end{equation}
and similarly for $\M{B}$.  By designating the null vector of
$\M{A}$ as $\V{u}$, we require that,
\begin{equation}
\M{M}\M{D}_y\M{N}\V{u} = \V{0}.
\end{equation}
Clearly, if we let
\begin{equation} \label{eqn:uDef}
\V{u} = \M{N}^{-1}\V{1},
\end{equation}
then we have,
\begin{eqnarray}
\M{A}\V{u} & = & \M{M}\M{D}_y\M{N}\M{N}^{-1}\V{1} \nonumber \\
& = & \M{M}\M{D}_y \V{1} \nonumber \\
& = & \V{0}
\end{eqnarray}
and thus $\V{u}$ as per Equation (\ref{eqn:uDef}) is the null
vector of the modified differential operator.
In the case where $\M{N}$ is orthonormal, as with spectral
regularization,
\begin{equation}
\V{u} = \MT{N}\V{1},
\end{equation}
and $\V{u}$ is the spectrum of a constant function. By the
identical derivation, we also have the null vector of the modified
differential operator $\M{B}$ as, $\M{B}\V{v} = \V{0}$.
The two null vectors define the null space of the Sylvester
Operator~\cite{Stewart2}, i.e., for
\begin{equation}
S\left(\M{\Phi}\right) = \MT{A}\M{A}\M{\Phi} +
\M{\Phi}\MT{B}\M{B},
\end{equation}
we have the null surface, $\M{\Phi}_0 = \alpha \V{u}\VT{v}$, such
that,
\begin{equation}
S\left( \M{\Phi}_0 \right) = \M{0},
\end{equation}
for arbitrary $\alpha$. The constant $\alpha$ is the effective
constant of integration. \\
 Thus, knowing the null surface of the Sylvester Operator, we
propose the following algorithm which essentially removes this
degree of freedom from the solution of the Sylvester Equation.
Introducing the Householder reflections~\cite{Golub},
\begin{equation}
\M{P}_a = \M{I}_{m} - 2
\frac{\V{\tilde{u}}\VT{\tilde{u}}}{\VT{\tilde{u}}\V{\tilde{u}}}
\qquad \textrm{and} \qquad
\M{P}_b = \M{I}_{n} - 2
\frac{\V{\tilde{v}}\VT{\tilde{v}}}{\VT{\tilde{v}}\V{\tilde{v}}}
\end{equation}
with
\begin{equation}
\V{\tilde{u}} = \V{u} + \left\| \V{u} \right\|_2 \V{e}_1
\qquad \textrm{and} \qquad
\V{\tilde{v}} = \V{v} + \left\| \V{v} \right\|_2 \V{e}_1
\end{equation}
and appropriately sized coordinate vectors, $\V{e}_1$, we
transform $\M{\Phi}$ such that,
\begin{equation}
\M{\Phi} = \M{P}_a \M{\Psi} \MT{P}_b.
\end{equation}
Substituting this expression into the Sylvester Equation in
Equation (\ref{eqn:genNE}), yields,
\begin{equation}
\MT{A}\M{A}\M{P}_a \M{\Psi} \MT{P}_b + \M{P}_a \M{\Psi}
\MT{P}_b\MT{B}\M{B} - \MT{A}\M{F} - \M{G}\M{B} = \M{0},
\end{equation}
By pre-multiplying by $\MT{P}_a$ and post-multiplying by
$\M{P}_b$, we yield the Sylvester Equation,
\begin{equation} \label{eqn:sylEqnHR}
\MT{\hat{A}}\M{\hat{A}}\M{\Psi} + \M{\Psi}\MT{\hat{B}}\M{\hat{B}}
- \MT{\hat{A}}\M{\hat{F}} - \M{\hat{G}}\M{\hat{B}} = \M{0}
\end{equation}
with\footnote{The following relations are equivalent to
differentiating the Householder reflections with the modified
differential operators.}
\begin{equation}
\M{\hat{A}} = \M{A}\M{P}_a = \begin{bmatrix} \V{0} & \M{R}
\end{bmatrix}, \qquad
\M{\hat{B}} = \M{B}\M{P}_b = \begin{bmatrix} \V{0} & \M{S}
\end{bmatrix}
\end{equation}
and
\begin{equation}
\M{\hat{F}} = \M{F}\M{P}_b, \qquad \M{\hat{G}} = \MT{P}_a\M{G}.
\end{equation}
The above structure of the matrices $\M{\hat{A}}$ and
$\M{\hat{B}}$ arises since the Householder reflections have been
chosen such that,
\begin{equation}
\M{\hat{A}}\V{e}_1 = \V{0}
\qquad \textrm{and} \qquad
\M{\hat{B}}\V{e}_1 = \V{0}.
\end{equation}
It is important to note however, that the Householder matrices
$\M{P}_a$ and $\M{P}_b$ should not be formed explicitly, as doing
so increases the relevant work by an order of
magnitude~\cite[p.211]{Golub}; all relevant information is
contained in the vectors $\V{\tilde{u}}$ and $\V{\tilde{v}}$.
By the above arguments, the ``right hand side" of the Sylvester
Equation (\ref{eqn:sylEqnHR}) takes the form,
\begin{eqnarray}
\MT{\hat{A}}\M{\hat{F}} + \M{\hat{G}}\M{\hat{B}} & = &
\begin{bmatrix} \VT{0} \\ \MT{R} \end{bmatrix}
\begin{bmatrix} \V{\hat{f}}_1 & \M{\hat{F}}_2 \end{bmatrix}
+
\begin{bmatrix} \VT{\hat{g}}_1 \\ \M{\hat{G}}_2  \end{bmatrix}
\begin{bmatrix} \V{0} & \M{S} \end{bmatrix} \\
& = & \begin{bmatrix} 0 & \VT{\hat{g}}_1\M{S} \\
\MT{R}\V{\hat{f}}_1 & \MT{R}\M{\hat{F}}_2 + \M{\hat{G}}_2\M{S}
\end{bmatrix},
\end{eqnarray}
whereby the first column is partitioned from $\M{\hat{F}}$ and the
first row is partitioned from $\M{\hat{G}}$.  Thus, the Sylvester
Equation in (\ref{eqn:sylEqnHR}) partitions as follows,
\begin{eqnarray}
\begin{bmatrix} 0 & \VT{0} \\ \V{0} & \MT{R}\M{R} \end{bmatrix}
\begin{bmatrix} \psi_{00} & \VT{\psi}_{01} \\ \V{\psi}_{10} & \M{\Psi}_{11} \end{bmatrix}
& + &
\begin{bmatrix} \psi_{00} & \VT{\psi}_{01} \\ \V{\psi}_{10} & \M{\Psi}_{11} \end{bmatrix}
\begin{bmatrix} 0 & \VT{0} \\ \V{0} & \MT{S}\M{S} \end{bmatrix} \nonumber \\
& = & \begin{bmatrix} 0 & \VT{\hat{g}}_1\M{S} \\
\MT{R}\V{\hat{f}}_1 & \MT{R}\M{\hat{F}}_2 + \M{\hat{G}}_2\M{S}
\end{bmatrix}
\end{eqnarray}
which represents four equations; due to this partitioning the set
of equations reads,
\begin{eqnarray}
\VT{\psi}_{01}\MT{S}\M{S} & = & \VT{\hat{g}}_1\M{S} \\
\MT{R}\M{R}\V{\psi}_{10} & = & \MT{R}\V{\hat{f}}_1 \\
\MT{R}\M{R} \M{\Psi}_{11} + \M{\Psi}_{11} \MT{S}\M{S} & = &
\MT{R}\M{\hat{F}}_2 + \M{\hat{G}}_2\M{S}. \label{eqn:fullRankSE}
\end{eqnarray}
The first two equations are the normal equations of two simple
linear systems, i.e., they represent the least squares solution to
the over-determined systems,
\begin{eqnarray}
\VT{\psi}_{01}\MT{S} & = & \VT{\hat{g}}_1 \\
\M{R}\V{\psi}_{10} & = & \V{\hat{f}}_1
\end{eqnarray}
These equations should be solved directly using an appropriate
least squares method (i.e., without forming the normal equations
as above, cf.~\cite{Golub}).
The remaining equation, i.e. Equation (\ref{eqn:fullRankSE}), is a
full rank Sylvester Equation in $\M{\Psi}_{11}$, and can therefore
be solved using a standard algorithm (e.g.,
\cite{bartels1972,golub1979}). The expression for $\psi_{00}$, is
$0 \psi_{00} = 0$, and hence it can be set arbitrarily to zero; it
represents the effective constant of integration for the
reconstruction problem. The solution of the rank deficient
Sylvester Equation is therefore,
\begin{equation}
\M{\Phi} = \M{P}_a \M{\Psi} \MT{P}_b
\end{equation}
with
\begin{equation}
\M{\Psi} =
\begin{bmatrix} 0 & \V{\psi}_{01}^{*\textrm{T}} \\ \V{\psi}_{10}^{*} & \M{\Psi}_{11}^{*} \end{bmatrix}
\end{equation}
where $\V{\psi}_{01}^{*}$, $\V{\psi}_{10}^{*}$, and
$\M{\Psi}_{11}^{*}$ represent computed solutions.  By setting
$\psi_{00}$ to zero, an implicit constraint on the parameters is
imposed such that,
\begin{equation}
\VT{u}\M{\Phi}\V{v} = 0.
\end{equation}
For the simple least squares solution, this means that,
\begin{equation}
\VT{1}\M{Z}\V{1} = 0,
\end{equation}
that is, that the reconstructed surface is mean free.  Similarly,
for the weighted least squares solution, the reconstructed surface
satisfies,
\begin{equation}
\VT{1}\M{\Lambda}_{xy}^{-1}\M{Z}\M{\Lambda}_{yx}^{-1}\V{1} = 0,
\end{equation}
that is, its weighted mean is zero.
%
%
\begin{table*}[!tb]
\renewcommand{\arraystretch}{1.3}
\caption{Key to Sylvester Equations} \label{tab:sylEqnKey}
\centering
\begin{tabular}{|l||c|c|c|c|c|c|c|}
\hline Algorithm & $\M{A}$ & $\M{B}$ & $\M{F}$ & $\M{G}$ &
$\M{\Phi}$ & $\M{Z}$ &
 $\mathop{\mathrm{null}}\left( S\left(\M{\Phi}\right) \right)$ \\
\hline \hline
%
GLS & $\M{D}_y$ & $\M{D}_x$ & $\M{\hat{Z}}_y$ & $\M{\hat{Z}}_x$ & $\M{Z}$ & $\M{Z}$ & $\V{u}=\V{1}$, $\V{v}=\V{1}$  \\
\hline
%
Spectral & $\M{D}_y\M{B}_y$ & $\M{D}_x\M{B}_x$ &
$\M{\hat{Z}}_y\M{B}_x$ & $\MT{B}_y\M{\hat{Z}}_x$ &
 $\M{C}$ & $\M{Z}=\M{B}_y\M{C}\MT{B}_x$ & $\V{u}=\V{e}_1$, $\V{v}=\V{e}_1$ \\
\hline
%
\begin{minipage}[c]{0.17\columnwidth} Tikhonov (Standard) \end{minipage} & $\begin{bmatrix} \M{D}_y \\
\lambda \M{I}_m
\end{bmatrix}$ &
  $\begin{bmatrix} \M{D}_x \\ \lambda \M{I}_n \end{bmatrix}$ &
  $\begin{bmatrix} \M{\hat{Z}}_y \\ \lambda \M{Z}_0 \end{bmatrix}$ &
  $\begin{bmatrix} \M{\hat{Z}}_x & \lambda \M{Z}_0 \end{bmatrix}$ &
  $\M{Z}$ & $\M{Z}$ & $\emptyset$ \\
\hline
%
\begin{minipage}[c]{0.17\columnwidth} Tikhonov (Degree-$k$) \end{minipage} & $\begin{bmatrix} \M{D}_y \\ \mu \M{D}_y^k
\end{bmatrix}$ & $\begin{bmatrix} \M{D}_x \\ \lambda \M{D}_x^k
\end{bmatrix}$ & $\begin{bmatrix} \M{\hat{Z}}_y \\ \mu \M{D}_y^k \M{Z}_0
\end{bmatrix}$ & $\begin{bmatrix} \M{\hat{Z}}_x & \lambda \M{Z}_0
\left(\MT{D}_x\right)^k
\end{bmatrix}$ & $\M{Z}$ & \M{Z} & $\V{u}=\V{1}$, $\V{v}=\V{1}$ \\
\hline
Dirichlet & $\M{D}_y\M{P}_m$ & $\M{D}_x\M{P}_n$ & $\left(
\M{\hat{Z}}_y - \M{D}_y\M{Z}_\textrm{B} \right)\M{P}_n$ &
$\MT{P}_m\left( \M{\hat{Z}}_x - \M{Z}_\textrm{B}\MT{D}_x
\right)$ & $\M{Z}_{\textrm{I}}$ & $\M{Z} = \M{P}_m\M{Z}_{\textrm{I}}\MT{P}_n + \M{Z}_\textrm{B}$ & $\emptyset$ \\
\hline
 Weighted & $\M{\Lambda}_{yy}^{-\frac{1}{2}}\M{D}_y
\M{\Lambda}_{xy}^{\frac{1}{2}}$ &
$\M{\Lambda}_{xx}^{-\frac{1}{2}}\M{D}_x
\M{\Lambda}_{yx}^{\frac{1}{2}}$ &
$\M{\Lambda}_{yy}^{-\frac{1}{2}}\M{\hat{Z}}_y
\M{\Lambda}_{yx}^{-\frac{1}{2}}$ &
$\M{\Lambda}_{xy}^{-\frac{1}{2}}\M{\hat{Z}}_x
\M{\Lambda}_{xx}^{-\frac{1}{2}}$ & $\M{Z}_\textrm{w}$ & $\M{Z} =
\M{\Lambda}_{xy}^{-\frac{1}{2}}\M{Z}_\textrm{w}\M{\Lambda}_{yx}^{-\frac{1}{2}}$
&
$\V{u}=\M{\Lambda}_{xy}^{-\frac{1}{2}}\V{1}$, $\V{v}=\M{\Lambda}_{yx}^{-\frac{1}{2}}\V{1}$ \\
\hline
%
%
\end{tabular}
\end{table*}
%
%
\section{Numerical Testing} \label{sec:numTest}
%
To demonstrate the functionality and purpose of the new
algorithms, the following numerical tests are proposed:
\begin{enumerate}
\item Average computation time of the new algorithms in comparison
to state-of-the-art algorithms.
\item Monte-Carlo simulations for the reconstruction of a surface
from its discrete-sampled analytic gradient field to demonstrate
the functionality in the presence of various types of noise.
\item Reconstruction from real Photometric Stereo data to
demonstrate the functionality with respect to a real-world
problem.
\end{enumerate}

\subsection{Computation Time}
%
The algorithm computation times have been computed by applying
each algorithm to solve an identical problem one-hundred times,
and averaging the results.  The computation times for existing
algorithms and the newly proposed are shown in Table
\ref{tab:timing}, which is divided into three parts:
\begin{enumerate}
\item Existing methods: GLS~\cite{Harker2008c}, GLS via sparse
matrix methods such as LSQR, and the methods of Frankot and
Chellappa~\cite{frankot1988}, and Horn and Brooks~\cite{horn1986}.
\item Computation time for an $n\times n$ Singular Value
Decompostion:  This provides a familiar reference with which to
compare the new algorithms.
\item Newly proposed methods: Spectral reconstruction, Tikhonov
regularization (known and unknown regularization parameter,
$\lambda$), Dirichlet boundary conditions, and the weighted least
squares solution.  Note that since the new methods are all direct,
the computation times are completely independent of the input data
(e.g., whether the gradient is smooth or completely random,
contains outliers, etc.).
\end{enumerate}
The state-of-the-art methods with
regularization~\cite{karacali2003,horovitz2004,agrawal2006,ng2010,balzer2011,balzer2012}
all fall under the category of GLS via sparse matrix methods; all
methods use some form of large-scale solver, and hence the
algorithms can be no faster than the times for solving the GLS
problem via LSQR.  The times for GLS (Sparse LSQR) thus provide an
order of magnitude estimate for the state-of-the-art methods.
However, some of these methods,
specifically~\cite{agrawal2006,ng2010,balzer2011,balzer2012} use
more elaborate approaches such as spline surfaces, which make them
far more computationally intensive than the sparse LSQR alone; due
to memory requirements and computation time these methods are
simply not functional on a modern PC for the surface sizes
presented here.
\begin{table}[htb]
\renewcommand{\arraystretch}{1.3}
\caption{Algorithm Computation Time (Seconds)} \label{tab:timing}
\centering
\begin{tabular}{|l||c|c|c|}
\hline
 & Small & Medium & Large \\
%
Algorithm & $2^7\times 2^7$ & $2^9\times 2^9$ & $2^{10}\times 2^{10}$ \\
\hline \hline
GLS & 0.0433 & 1.9417 & 14.4622 \\
\hline
GLS (Sparse LSQR) & 0.4782 & 43.7778 & 338.0722 \\
\hline
Frankot-Chellappa & 0.0035 & 0.0670 & 0.3992 \\
\hline
Horn-Brooks & 0.9858 & 12.5509 & 56.2055 \\
\hline \hline $n\times n$ SVD & 0.0283 & 2.4949 & 20.4829 \\
\hline \hline
Spectral & 0.0107 & 0.3455 & 2.5286 \\
\hline
Dirichlet & 0.0333 & 1.4328 & 11.7184 \\
\hline
Tikhonov (known $\lambda$)& 0.0423 & 1.9561 & 16.5140 \\
\hline
Weighted & 0.0602 & 2.7328 & 20.2516 \\
\hline
Tikhonov (L-curve) & 0.0700 & 6.0875 & 48.1106 \\
\hline
\end{tabular}
\end{table}
Of the existing methods, the Frankot-Chellappa algorithm is
clearly the fastest.  However, this method, along with the
Horn-Brooks algorithm can be considered to be approximate methods.
This can be demonstrated by the fact that the Spectral Method
proposed here using the Fourier basis yields exact reconstruction,
whereby the Frankot-Chellappa algorithm cannot.  The results
exhibit a low-frequency bias, and generally the results are
peculiar and unusable for non-periodic data; its computational
efficiency is hence not advantageous in any way.  Similarly, the
Horn-Brooks method is generally non-convergent; the times in Table
\ref{tab:timing} represent the time for 1000 iterations, whereby
several thousand are required to obtain a reasonable result (ca.
8000~\cite{durou2007}). The starkest contrast in Table
\ref{tab:timing} is between the GLS algorithm via the Sylvester
Equation, i.e., Equation (\ref{eqn:normalEq}), and solving the
exact same problem using sparse matrix methods, i.e., Equation
(\ref{eqn:normalEqVec}). The Sylvester Equation method
reconstructs the small surface in under 50ms, and the large
$1024\times 1024$ surface in a reasonable amount of time. The
times for the sparse solver are clearly an order of magnitude
larger than the standard GLS algorithm.  Thus, the
state-of-the-art methods with
regularization~\cite{karacali2003,horovitz2004,agrawal2006,ng2010}
are already at a large handicap with respect to the Sylvester
Equation methods presented here.  For example, Ng et
al.~\cite{ng2010} reported a time of $11856.20$s (about
$3\frac{1}{2}$ hours) to reconstruct a $240\times 314$ surface.\\
The lower portion of Table \ref{tab:timing} shows the average
computation times for the new algorithms; they are ordered in
terms of computational demand.  For example, the Spectral
Regularization method and Dirichlet Method requires only the
solution of a Sylvester Equation; the remaining methods require
the additional algorithm to account for rank-deficiency of the
Sylvester Equation; the weighted least squares algorithm requires
the computation of matrix square roots; and finally, the L-Curve
method of Tikhonov regularization requires additional functional
evaluations depending on the number of points desired on the
L-Curve.  Note that the majority of the algorithms are faster than
the $m\times n$ SVD computation because the solution of the
Sylvester Equation does not require full diagonalization of the
coefficient matrices; the exceptions are the Tikhonov L-Curve
algorithm which uses two SVDs, and the weighted least squares
which requires four matrix square roots. Clearly the most
efficient algorithm is the Spectral Regularization for reasons
discussed in Section \ref{sec:specReg}, whereby half of the basis
functions were truncated. Since the high frequency components of
sinusoidal functions or polynomials can most often be considered
noise, this algorithm appears to be the most advantageous.  The
Dirichlet and Tikhonov with known $\lambda$ algorithms yield times
which are on par with the standard GLS solution.  The matrix
square root computation adds a tangible overhead to the weighted
least squares algorithm.  The most computationally intensive new
algorithm is the Tikhonov regularization method with the L-Curve
to determine $\lambda$.  It required just over 48s to compute ten
points on the L-Curve, find the optimal $\lambda$ and reconstruct
the surface; clearly the computation is dominated by the two SVD
computations, rather than the norm computations. This is more than
reasonable given the sheer difficulty in determining
regularization parameters~\cite{Engl,GolubMeurant}. In any case,
the algorithm is three orders of magnitude faster (i.e., 1000
times faster for the $512\times 512$ surface) than the method of
Ng et al. for their $240\times 314$ surface.
\subsection{LS Properties via Monte-Carlo}

\begin{figure}[ht]
\centering
\subfigure[Test Surface $\M{Z}$]{
\includegraphics[clip=true]{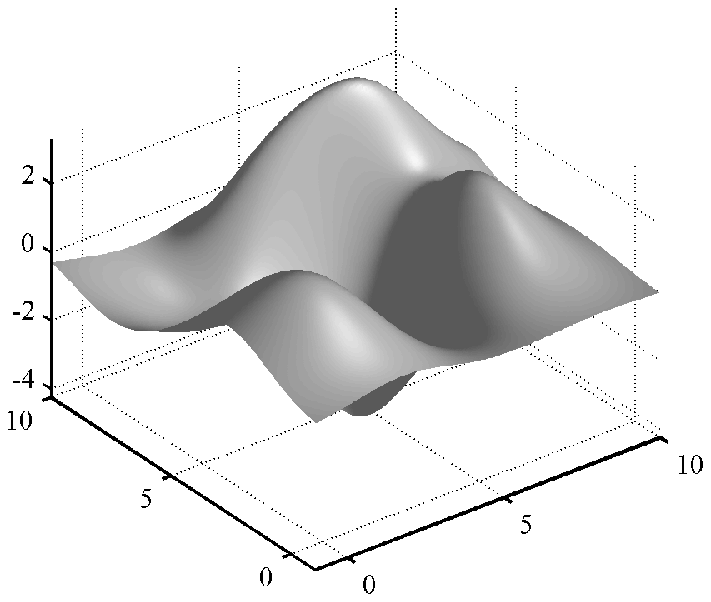}
\label{fig:mcSurfsub} }
\subfigure[$\M{Z}_x$]{
\includegraphics[clip=true,width=0.45\columnwidth]{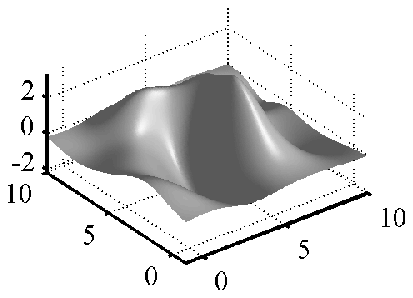}
\label{fig:mcGradX} }
\subfigure[$\M{Z}_y$]{
\includegraphics[clip=true,width=0.45\columnwidth]{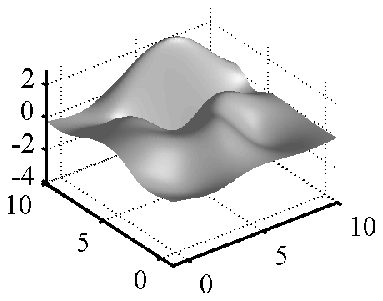}
\label{fig:mcGradY} }
\caption[-]{Ground truth for the Monte-Carlo simulation
\subref{fig:mcSurfsub} The $150\times 150$ test surface $\M{Z}$
and \subref{fig:mcGradX}, \subref{fig:mcGradY} its analytic
gradient field $\M{Z}_x$ and $\M{Z}_y$.} \label{fig:mcSurf}
\end{figure}

To demonstrate the functionality of the algorithms proposed here,
Monte-Carlo testing has been undertaken with various forms of
synthetic noise.  The test surface is an analytic surface which is
the sum of anisotropic Gaussian probability density functions of
the form,
\begin{equation}
z(x,y) = \sum_{k=1}^n A_k \exp \left( -\frac{1}{2} \VT{x}
\M{\Lambda}_k^{-1} \V{x} \right),
\end{equation}
which is similar to MATLAB$^{\circledR}$'s ``peaks" test function.
The test function and its gradient field are shown in Figure
\ref{fig:mcSurf}.
The motivation for using such a test function is twofold: the
reconstruction can be performed on the analytic derivatives
evaluated on a discrete grid; the surface is non-polynomial, which
means that only high order derivative approximations will give
accurate results.  In so doing, the reconstructed surfaces can be
compared in terms of relative error with respect to the exact
surface, and hence we obtain a quantitative measure of the quality
of the surface reconstruction.  It should be noted that this is
\textit{never} done in the literature, mainly due to the fact that
most published algorithms have some form of systematic error which
would skew such results.  Therefore, in the literature, it is
common to find only subjective reconstruction results, leaving the
reader to ``eyeball" the quality of the results. The Monte-Carlo
experiment, however, is exceptionally illuminating with respect to
the insight it gives into the functionality of each algorithm.  It
is therefore quite peculiar that researchers have yet to perform
such experiments. Indeed, with such excessive computation time
required for the vectorized solutions, a Monte-Carlo type
simulation is all but precluded with previous methods. The
Monte-Carlo simulations here, originally proposed
in~\cite{Harker2008c,harker2011}, thus represent the first ever
attempt to benchmark
solutions to the reconstruction of a surface from its gradient.\\
  As for the various forms of regularization proposed in this
paper, their strengths and weaknesses can be characterized in
terms of the type of noise present in the measurement.  Hence
Monte-Carlo testing has been performed with three types of noise:
\begin{enumerate}
\item The gradient field corrupted by i.i.d. Gaussian noise, in
which case the GLS algorithm is the ``Gold Standard" benchmark
solution.
\item The gradient field corrupted by heteroscedastic Gaussian
noise, in which case the weighted least squares solution is
optimal.
\item The gradient field corrupted by gross-outliers.  The
outliers are placed randomly throughout the gradient field such
that a percentage value of the gradient field is corrupted.
Outlier values are set to the maximum value of the respective
gradient component to mimic the saturated pixels of an image.
\end{enumerate}
\subsubsection{I.I.D.\ Gaussian Noise}
\begin{figure}[ht]
\centering
\subfigure[$\M{Z}_x$]{
\includegraphics[clip=true]{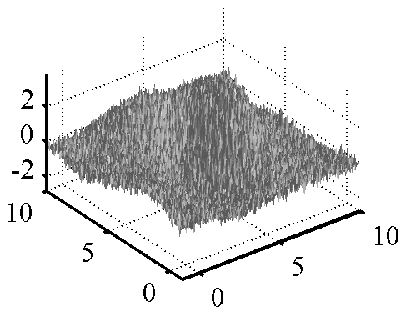}
\label{fig:mcGradXN} }
\subfigure[$\M{Z}_y$]{
\includegraphics[clip=true]{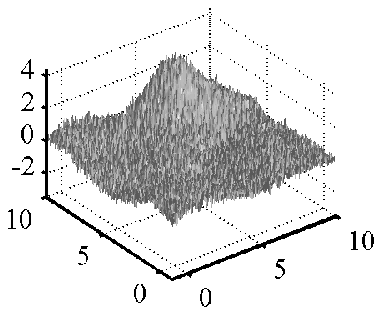}
\label{fig:mcGradYN} }
\caption[-]{Noisy gradient field for the Monte-Carlo simulation,
corrupted by i.i.d.\ Gaussian noise with a standard deviation of
$10\%$ of the gradient amplitude.} \label{fig:mcSurfNoisy}
\end{figure}

\begin{figure}[htb]
\centering
\includegraphics[clip=true,width=\columnwidth]{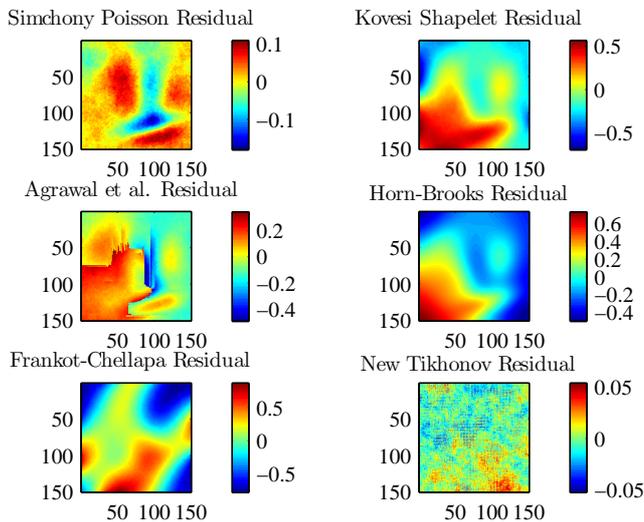}
\caption{Residuals of various methods: Existing methods exhibit
large systematic errors in their residuals.  Only the newly
proposed methods, such as the Tikhonov regularized solution,
exhibit purely stochastic residuals, which one would expect from a
least squares solution.} \label{fig:residuals}
\end{figure}

\begin{figure}[htb]
\centering
\includegraphics[clip=true,width=\columnwidth]{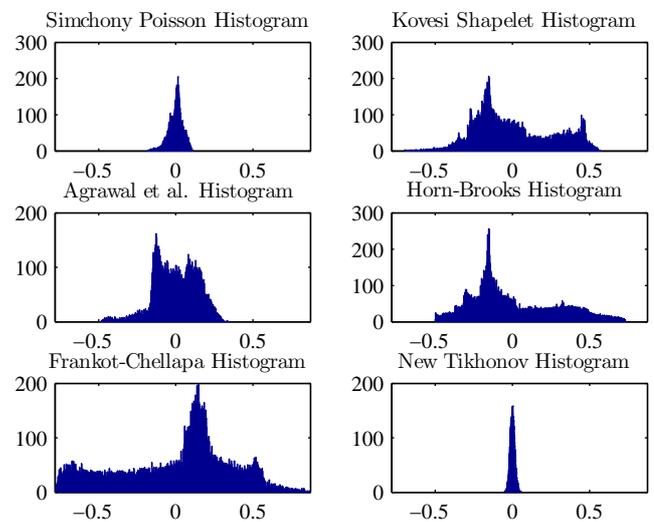}
\caption{Histograms of the residuals of various methods: only the
newly proposed method produces residuals which are themselves
Gaussian; note that they are also much smaller in magnitude.}
\label{fig:histograms}
\end{figure}

\begin{figure}[htb]
\centering
\includegraphics[clip=true,width=\columnwidth]{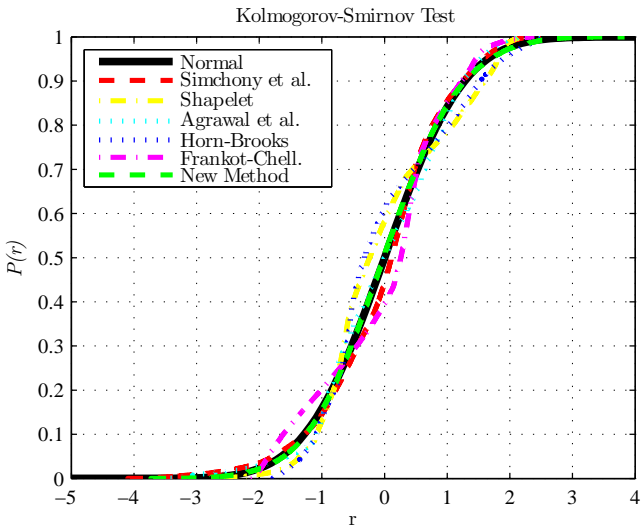}
\caption{Results from the Kolmogorov-Smirnov test. Clearly, only
the new method has Gaussian residuals.  This demonstrates
definitively that previous methods do not solve the least squares
reconstruction problem.} \label{fig:ksTest}
\end{figure}

\begin{figure}[htb]
\centering
\includegraphics[clip=true]{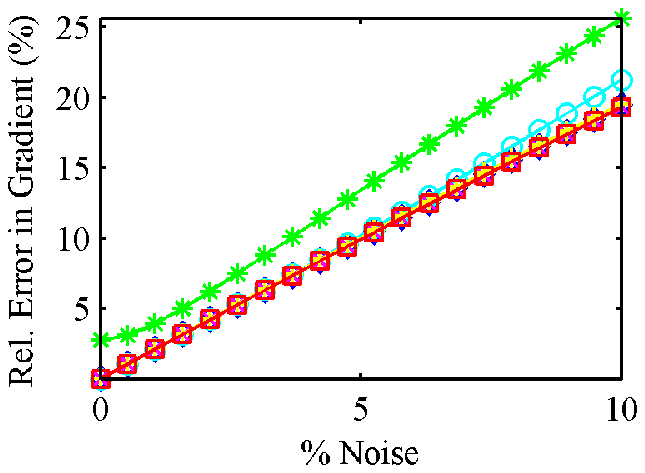}
\includegraphics[clip=true]{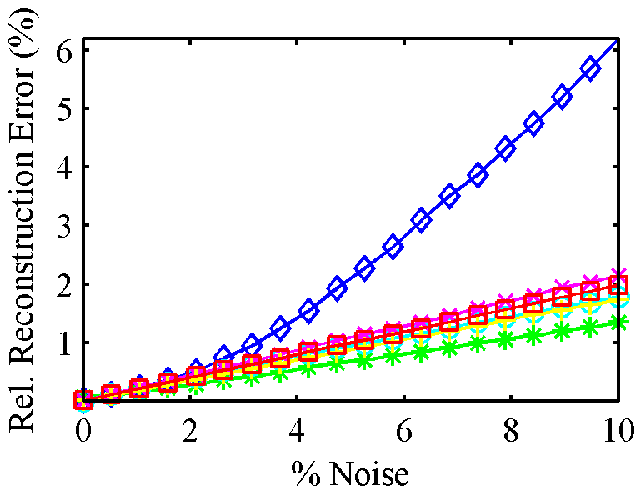}
\caption{Reconstruction subject to i.i.d.\ Gaussian noise. (TOP)
Relative cost function residual (BOT) Relative Reconstruction
Error. Legend: GLS ($\Box$), Spectral-Cosine ($\ast$), Tikhonov
Standard Form ($\diamond$), Tikhonov Degree-2 ($\circ$), Weighted
Least Squares ($\times$), Dirichlet Boundary Conditions ($+$).}
\label{fig:mcIID}
\end{figure}

\begin{figure*}[ht]
\centering
\subfigure[Global Least Squares (GLS)]{
\includegraphics[clip=true]{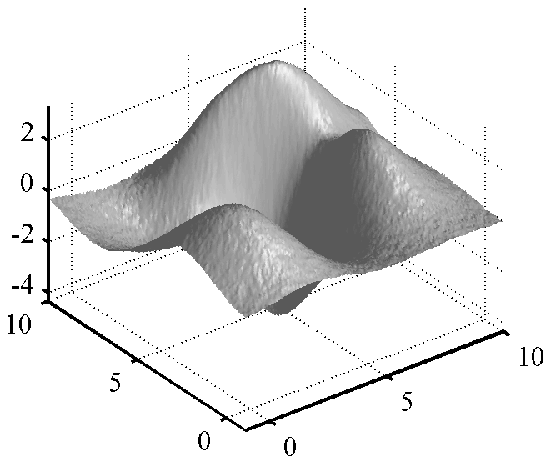}
\label{fig:mcIIDGLS} }
\subfigure[Spectral-Cosine]{
\includegraphics[clip=true]{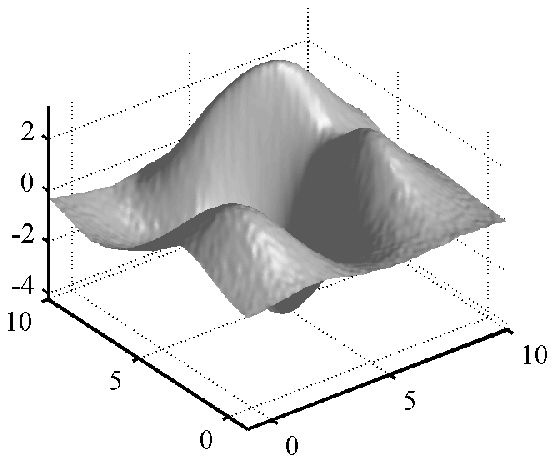}
\label{fig:mcIIDSC} }
\subfigure[Tikhonov Degree-2]{
\includegraphics[clip=true]{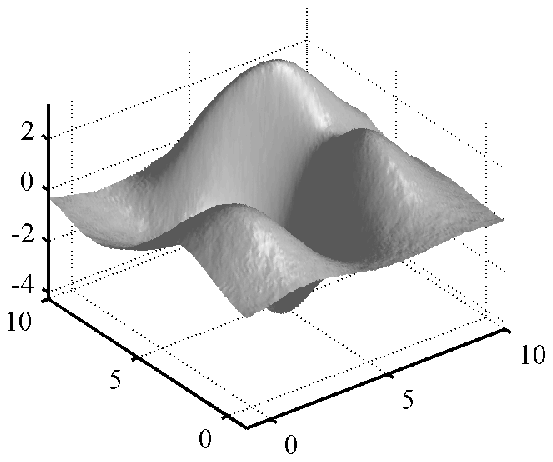}
\label{fig:mcIIDT2} }
\caption[-]{Reconstructions at maximum noise for i.i.d.\ Gaussian
noise \subref{fig:mcIIDGLS} the GLS solution \subref{fig:mcIIDSC}
the spectral-cosine reconstruction with half of the basis
functions truncated \subref{fig:mcIIDT2} Tikhonov reconstruction
with degree-2 regularization term.} \label{fig:mcIIDSurfs}
\end{figure*}

Given the existing methods in the literature, the reconstruction
of a surface from a gradient corrupted by i.i.d.\ Gaussian noise
has been, until now, a notoriously difficult problem.
Figure~\ref{fig:mcSurfNoisy} shows the gradient of the surface in
Figure~\ref{fig:mcSurf} when corrupted with Gaussian noise with a
standard deviation $10\%$ of the gradient amplitude.
Figure~\ref{fig:residuals} shows the reconstruction residuals of
various existing
methods~\cite{simchony1990,kovesi2005,agrawal2006,horn1986,frankot1988}
as compared to the new Tikhonov solution.  All the existing
methods exhibit substantial systematic error in their residuals;
in comparison, the new Tikhonov solution has a residual matrix
which is purely stochastic -- a typical feature of a least squares
solution proper.  In Figure \ref{fig:histograms}, the histograms
of these residuals are plotted.  The existing methods exhibit
highly irregular distributions due to the systematic errors in
their computation.  The residuals of the Tikhonov solution
presented here are firstly Gaussian, and secondly significantly
smaller than those of the state-of-the-art solutions.  These are
the results one would expect, statistically speaking, from an
appropriate global least squares solution.  Specifically, in
Figure \ref{fig:histograms} the method of Simchony shows the
results of using poor or incorrect derivative formulas; the method
of Kovesi (akin to the Frankot-Chellappa method) shows the results
of using inappropriate basis functions; and the method of Agrawal
et al. shows the general inappropriateness of path integration
methods when noise is present in the data.  The other methods
demonstrate similarly biased results. Similar skewed non-Gaussian
residuals are obtained with the FEM type methods such as the
method of Balzer and M\"{o}rwald~\cite{balzer2012}.
  The final ``nail-in-the-coffin" of previous methods is
demonstrated by means of a Kolomogorov-Smirnov statistical test.
Figure \ref{fig:ksTest} shows the normalized distributions of the
residuals of various methods.  If the methods yield Gaussian
residuals, then their density functions should be close to that of
a normal distribution.  In this case, the residuals of the new
method are almost indistinguishable from the normal distribution.
In contrast, none of the previous methods come close to a normal
distribution.  This demonstrates indisputably, that the global
least squares surface reconstruction from gradients problem has
been, until now, an unsolved problem.  In light of the fact that
the existing methods have extremely poor noise properties, they
have not been included in the subsequent Monte-Carlo tests;
namely, their residuals are typically an order of magnitude
larger, and hence would obscure any graphical method of
comparison.
\\
When the gradient field is corrupted by i.i.d.\ Gaussian noise,
the maximum likelihood reconstruction is provided by the least
squares solution.  That is to say, the ``Gold-Standard" in this
case is the GLS solution, in that is attains the lowest possible
bound of the cost function; it is hence the benchmark solution in
the presence of i.i.d.\ Gaussian noise.
  The relative value of the cost function
attained and the relative reconstruction error are shown for the
Monte-Carlo simulation in Figure \ref{fig:mcIID}.  The lower bound
of the cost function is attained by the least squares solution.
Notable is the cost function residual of the spectral methods;
these have the largest cost function residual, however, they
provide the best reconstruction residual.  This is due to the fact
that the high frequency components have been eliminated, and in
this case they correspond to only noise.  In contrast, with the
least squares solution, this high frequency noise is to some
degree integrable.  Also of note is the reconstruction with
Standard Tikhonov regularization.  The reconstruction is not as
good because systematic low-frequency errors are suppressed, of
which there are none in this case.  The degree-2 Tikhonov provides
better results as it suppresses the high frequency components. The
reconstructed surfaces for the GLS, Spectral-Cosine, and degree-2
Tikhonov methods, all at maximum noise level are shown in Figure
\ref{fig:mcIIDSurfs}.  The GLS solution exhibits some texture due
to the high level of Gaussian noise.  The Spectral-Cosine and
degree-2 Tikhonov successfully smooth these high frequency
components.
%
%
\subsubsection{Heteroscedastic Gaussian Noise}
\begin{figure}[htb]
\centering
\includegraphics[clip=true]{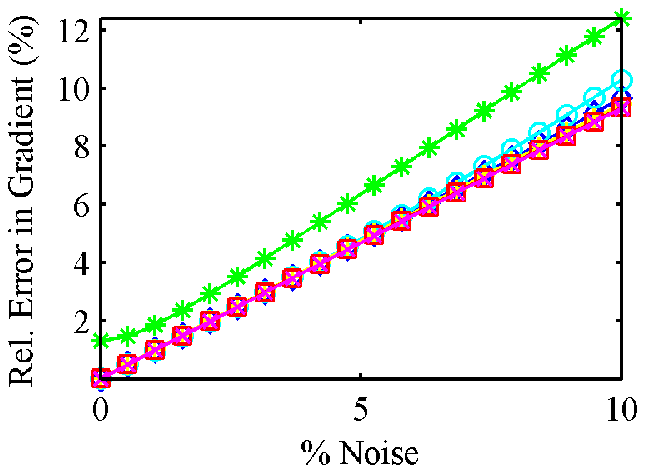}
\includegraphics[clip=true]{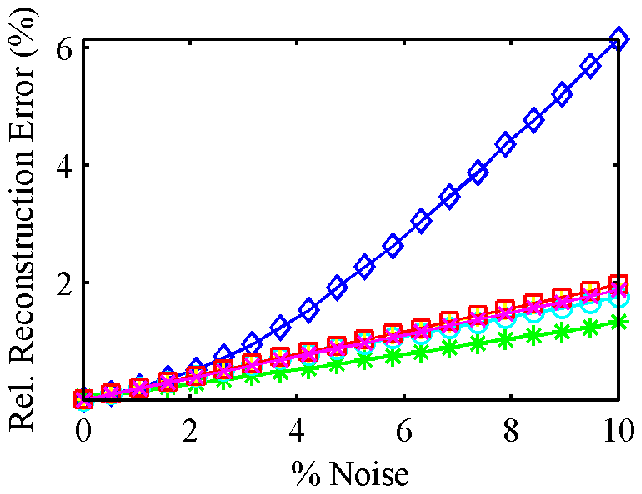}
\caption{Reconstruction subject to covariant Gaussian noise, in
which case weighted least squares is optimal. (TOP) Relative cost
function residual (BOT) Relative Reconstruction Error. Legend: GLS
($\Box$), Spectral-Cosine ($\ast$), Tikhonov Standard Form
($\diamond$), Tikhonov Degree-2 ($\circ$), Weighted Least Squares
($\times$), Dirichlet Boundary Conditions ($+$).}
\label{fig:mcCOV}
\end{figure}

\begin{figure*}[ht]
\centering
\subfigure[Weighted Least Squares (WLS)]{
\includegraphics[clip=true]{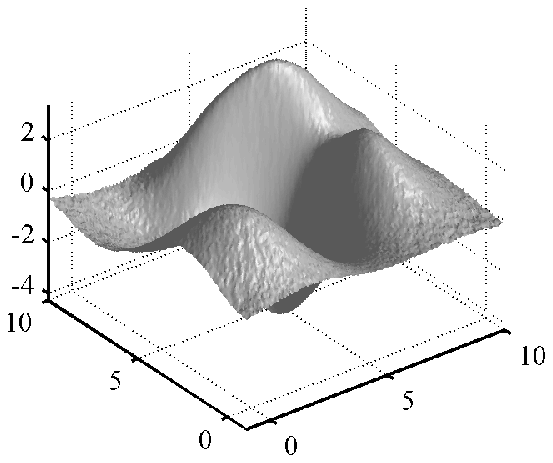}
\label{fig:mcCOVWLS} }
\subfigure[Spectral-Cosine]{
\includegraphics[clip=true]{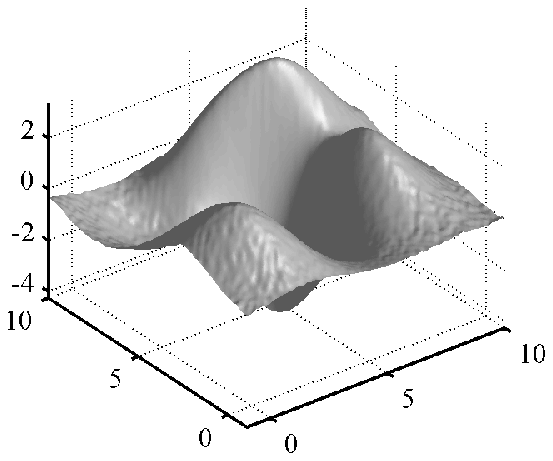}
\label{fig:mcCOVSC} }
\subfigure[Tikhonov Degree-1]{
\includegraphics[clip=true]{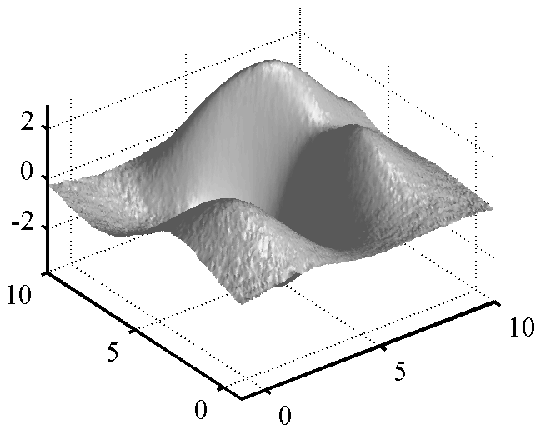}
\label{fig:mcCOVT2} }
\caption[-]{Reconstructions at maximum noise for covariant
Gaussian noise \subref{fig:mcCOVWLS} the WLS solution
\subref{fig:mcCOVSC} spectral cosine regularization with half the
basis functions truncated \subref{fig:mcCOVT2} reconstruction with
degree-1 Tikhonov regularization.} \label{fig:mcCOVSurfs}
\end{figure*}

If the noise in the gradient field is anisotropic, then the
maximum likelihood solution is given by the weighted least squares
solution.  For the Monte-Carlo test, the gradient field was
corrupted by a radially symmetric noise distribution with
increasing noise amplitude towards the image edges; this mimics
the error induced in photometric stereo by making the orthographic
projection assumption.  Results of the Monte-Carlo simulation are
shown in Figure \ref{fig:mcCOV}.  Clearly the weighted least
squares solution defines the lower bound of the cost function.
Again, the standard form Tikhonov regularization provides
relatively poor reconstruction since there is no systematic error
present.  Similarly to the i.i.d.\ case, the spectral methods have
the largest cost function residual, but due to their low-pass
functionality, again provide the best reconstruction.
In Figure \ref{fig:mcCOVSurfs}, the reconstructions at maximum
noise are shown for the weighted least squares, the spectral
reconstruction, and Tikhonov degree-1.  The Tikhonov degree-1
solution has the effect of suppressing the undulations of the
surface, which can have a similar effect as the WLS solution.
Clearly in the case of anisotropic noise, a weighted least squares
with spectral regularization would provide an optimal solution.
This can be accomplished with weighted basis
functions~\cite{oleary2010a}.
%
%
\subsubsection{Gross Outliers}
\begin{figure}[htb]
\centering
\includegraphics[clip=true]{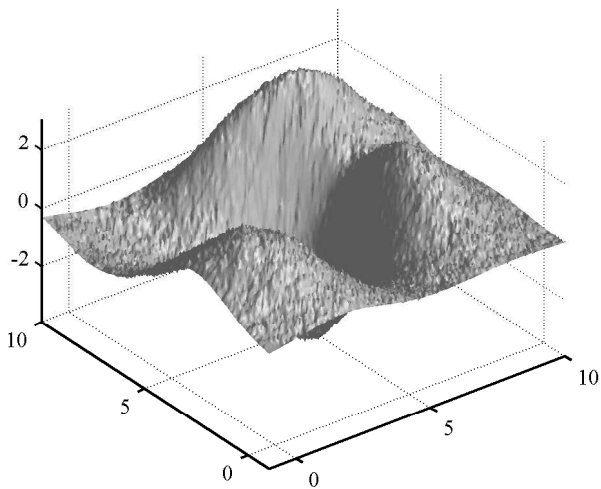}
\caption{Dirichlet reconstruction with $10\%$ outliers.  Given
this is an outlandish amount of noise, the Dirichlet
reconstruction is shown to be exceptionally robust to outliers.}
\label{fig:mcOUTdirSurf}
\end{figure}

\begin{figure}[htb]
\centering
\includegraphics[clip=true]{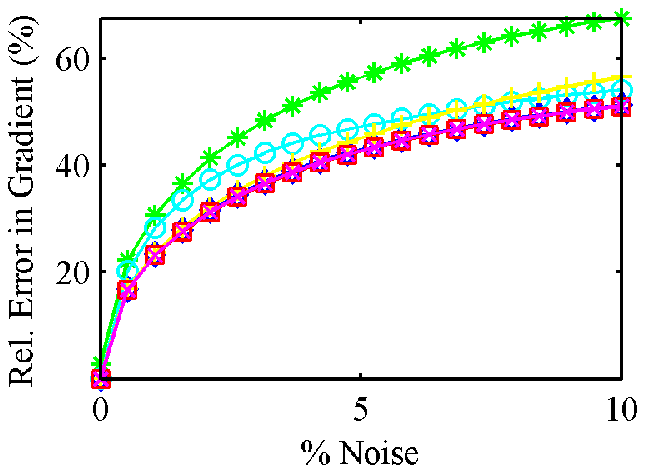}
\includegraphics[clip=true]{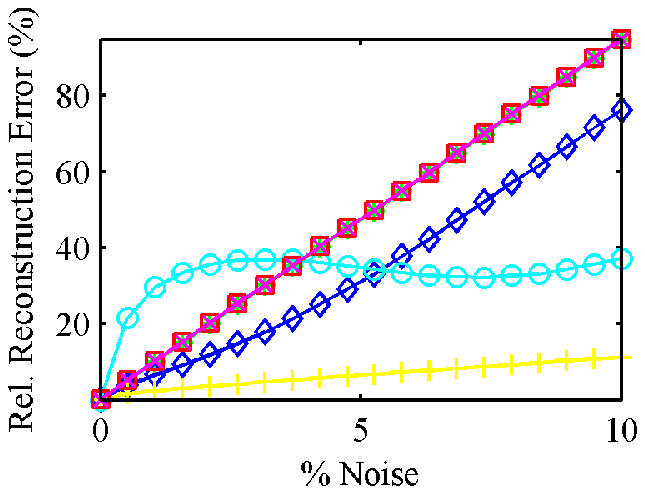}
\caption{Reconstruction subject gross outliers (saturated gradient
pixels). (TOP) Relative cost function residual (BOT) Relative
Reconstruction Error. Legend: GLS ($\Box$), Spectral-Cosine
($\ast$), Tikhonov Standard Form ($\diamond$), Tikhonov Degree-2
($\circ$), Weighted Least Squares ($\times$), Dirichlet Boundary
Conditions ($+$).} \label{fig:mcOUT}
\end{figure}

\begin{figure*}[ht]
\centering
\subfigure[Spectral-Cosine Low-Pass]{
\includegraphics[clip=true,width=0.31\textwidth]{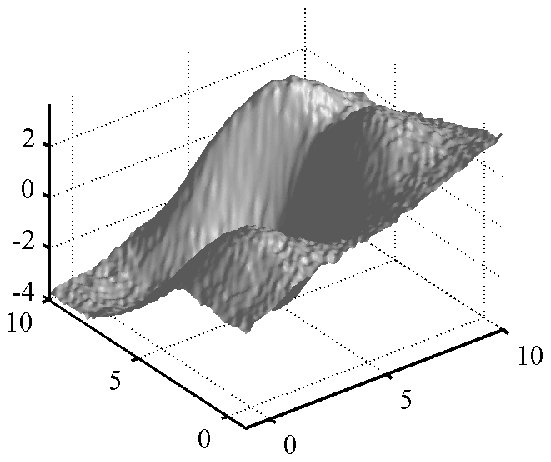}
\label{fig:mcOUTSLP} }
\subfigure[Spectral-Gram Band-Pass]{
\includegraphics[clip=true,width=0.31\textwidth]{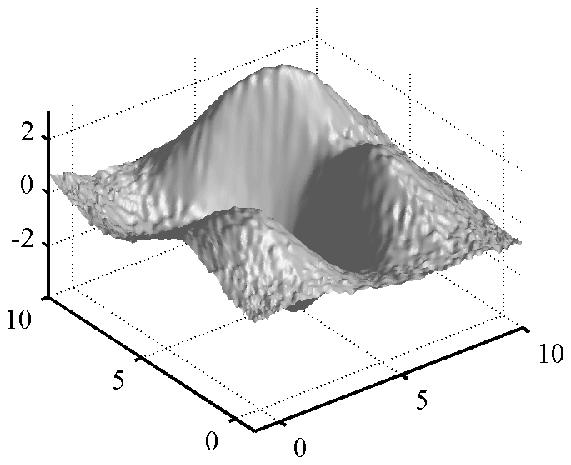}
\label{fig:mcOUTSBP} }
\subfigure[Tikhonov Standard Form]{
\includegraphics[clip=true,width=0.31\textwidth]{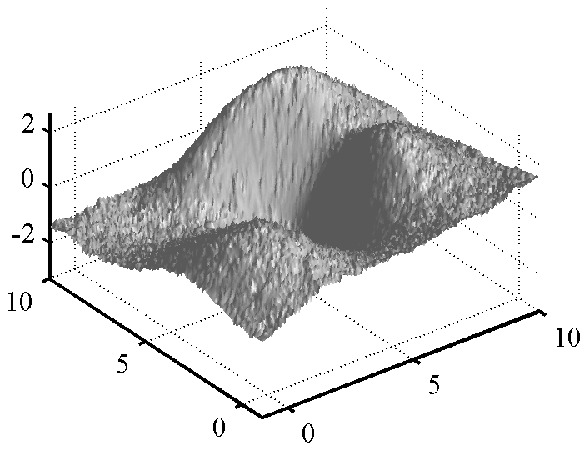}
\label{fig:mcOUTT0} }
\caption[-]{Reconstructions at maximum noise for gross outliers
\subref{fig:mcOUTSLP} Spectral-Cosine with half the basis
functions truncated \subref{fig:mcOUTSBP} Spectral-Gram with half
the basis functions truncated as well as the linear components
\subref{fig:mcOUTT0} Tikhonov Standard form with $\lambda$
determined by the L-Curve.} \label{fig:mcOUTSurfs}
\end{figure*}

To demonstrate the functionality of the algorithms in the presence
of outliers, a Monte-Carlo test was performed based on percentage
of outliers.  That is, for a given percentage of outliers, random
pixels in the gradient were set to the maximum amplitude of the
gradient component to simulate saturated image pixels; this
simulates what can transpire in real photometric stereo when
specular reflection occurs on the object surface.  In this case,
Tikhonov regularization in standard form (degree-0) algorithm is
optimal, since the outliers create a low frequency systematic bias
in the solution.  Results of the Monte-Carlo simulation are shown
in Figure (\ref{fig:mcOUT}), which shows the cost function
residual and reconstruction error as a function of percentage of
outliers.  Clearly the cost function residuals do not follow the
linear trend which is typical of a least squares solution subject
to Gaussian noise.  In the case of the reconstruction error, by
far the best reconstruction is provided by the Dirichlet boundary
conditions; clearly, if the value of the surface at the boundary
is known, the reconstruction can be extremely robust to outliers.
The reconstruction with Dirichlet boundary conditions is shown in
Figure \ref{fig:mcOUTdirSurf}.
In this case the low-pass spectral reconstruction is rather poor;
however, a band-pass spectral reconstruction can be used to remove
low-frequency components largely due to the outliers.  The
standard form Tikhonov regularization successfully suppresses much
of the bias due to the presence of outliers.  The reconstruction
results at maximum noise (percent outliers) are shown in Figure
\ref{fig:mcOUTSurfs} for low-pass spectral reconstruction,
band-pass Spectral-Gram reconstruction, and degree-0 Tikhonov
regularization. The large bias of the low-pass spectral
reconstruction is evident; the ``saturated" outliers induce a
large DC component into the gradient, which integrates to a ramp
function.  However both the band-pass reconstruction (removing all
linear polynomial components) and the Tikhonov regularization
successfully remove systematic bias from the solution.
\subsection{Real Photometric Stereo}
\begin{figure}[htb]
\centering
\includegraphics[clip=true]{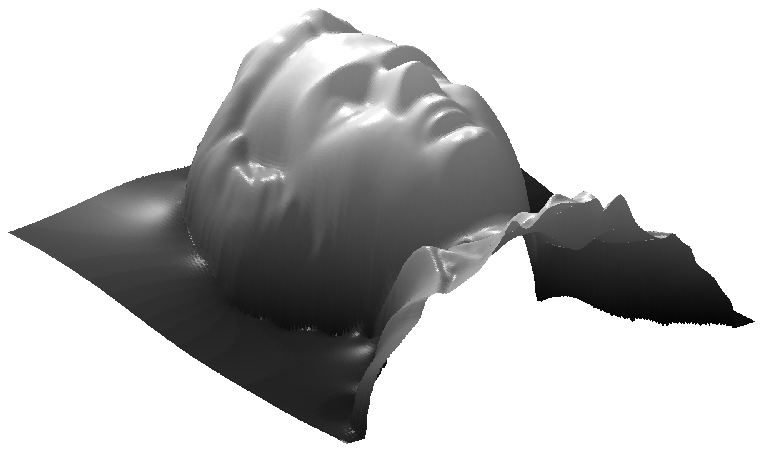}
\includegraphics[clip=true]{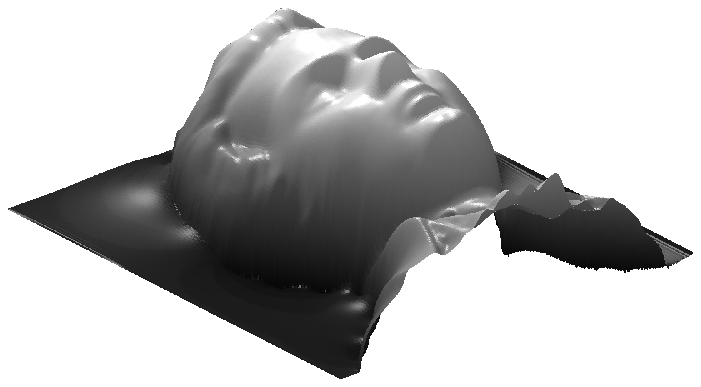}
\caption{The ``Mozart" data set reconstructed with (TOP) the GLS
solution and (BELOW) reconstruction with homogeneous Dirichlet
boundary conditions on three sides.} \label{fig:mozart}
\end{figure}

\begin{figure}[htb]
\centering
\includegraphics{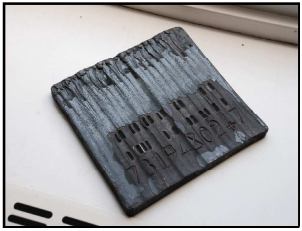}
\caption{Surface whose gradient was measured via Photometric
Stereo.} \label{fig:measuredSurf}
\end{figure}

\begin{figure*}[ht]
\centering
\subfigure[Global Least Squares (GLS)]{
\includegraphics[clip=true]{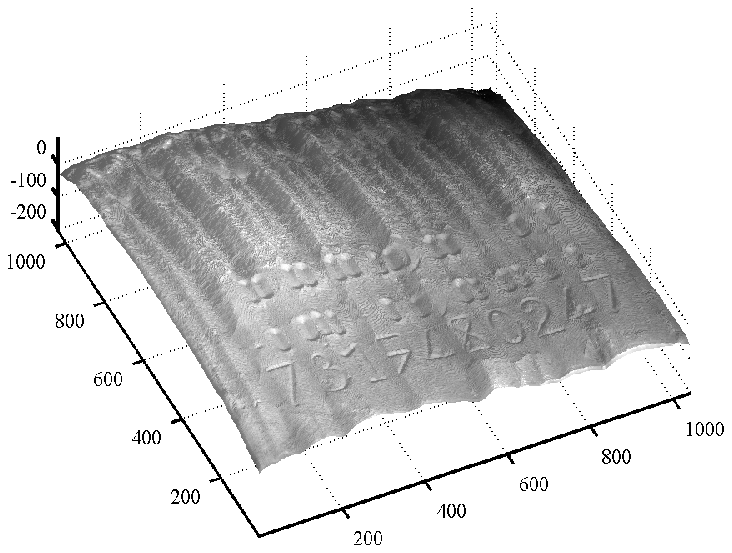}
\label{fig:realGLS} }
\subfigure[Spectral-Gram Band Pass]{
\includegraphics[clip=true]{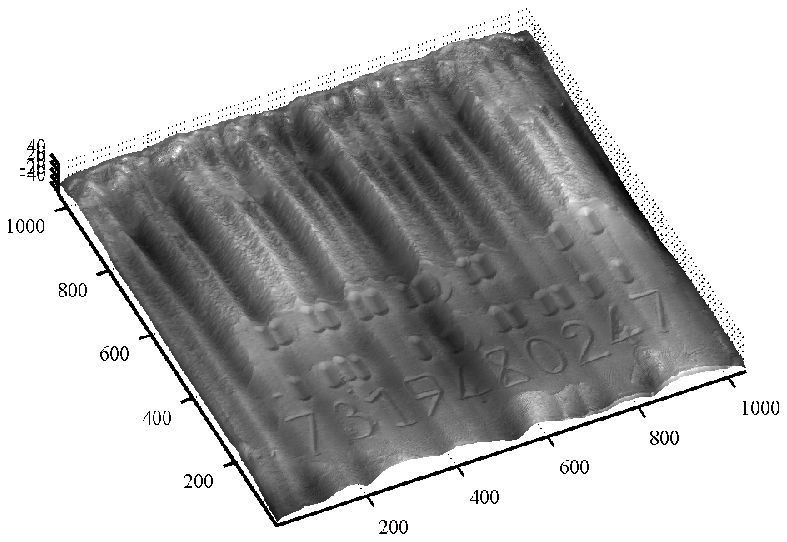}
\label{fig:realSGBP} }
\subfigure[Tikhonov Standard Form $\lambda^*$ from L-Curve]{
\includegraphics[clip=true]{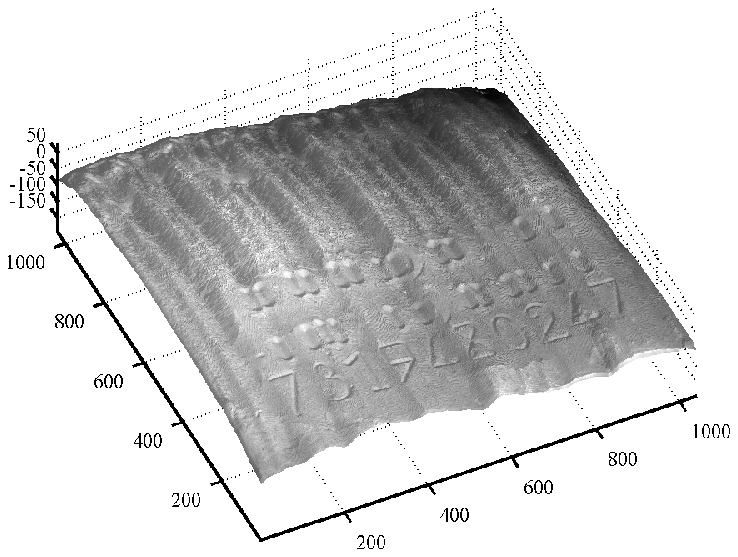}
\label{fig:realTikLam} }
\subfigure[Tikhonov Standard Form $\lambda=3\lambda^*$]{
\includegraphics[clip=true]{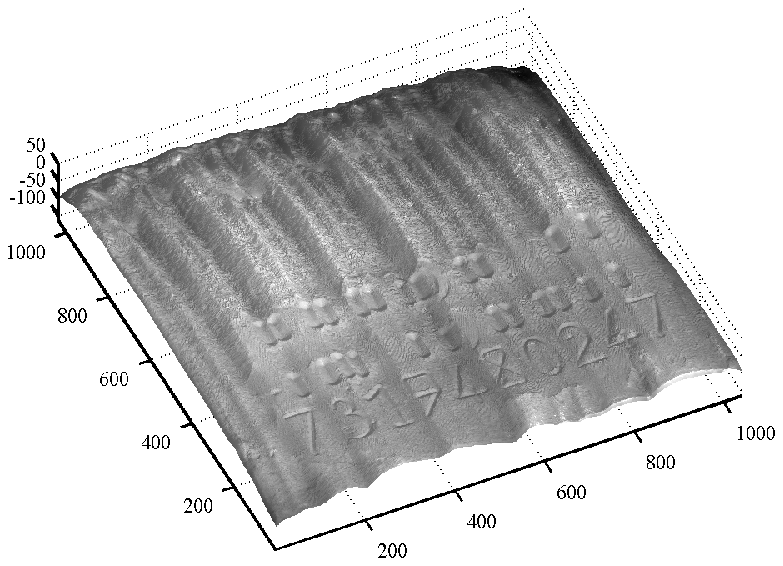}
\label{fig:realTik3Lam} }
\subfigure[Weighted Least Squares (WLS)]{
\includegraphics[clip=true]{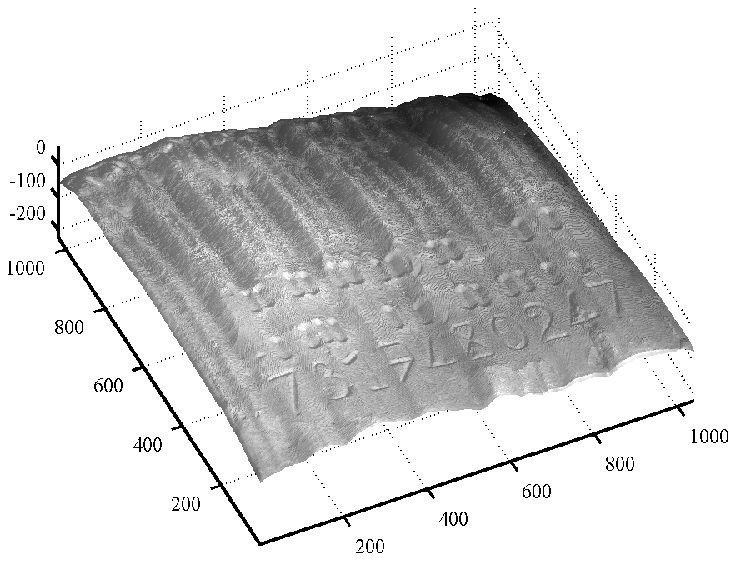}
\label{fig:realWLS} }
\subfigure[Homogeneous Dirichlet Boundary Conditions]{
\includegraphics[clip=true]{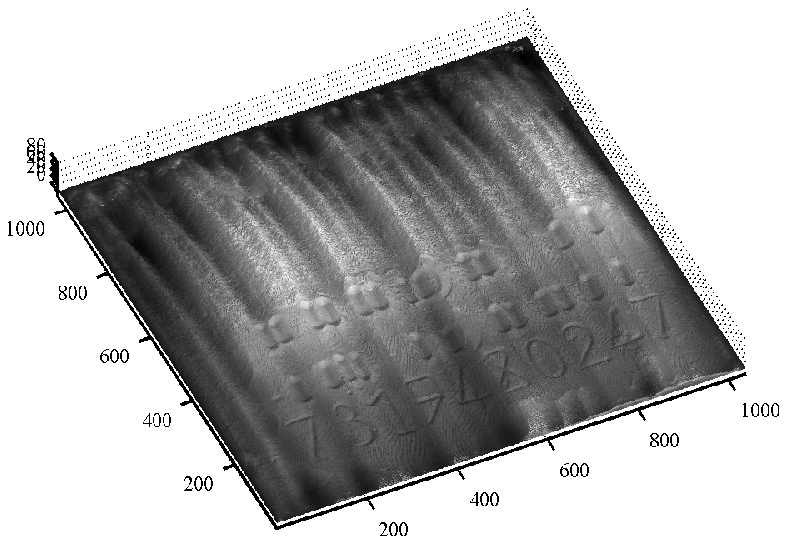}
\label{fig:realDIR} }
\caption[-]{Reconstructions of a real data set obtained via
Photometric Stereo \subref{fig:realGLS} GLS solution
\subref{fig:realSGBP} Spectral-Gram with all cubic terms removed,
as well as the high-frequency half of the basis functions
\subref{fig:realTikLam} Tikhonov Standard Form with $\lambda^*$
determined by the L-Curve \subref{fig:realTik3Lam} Tikhonov
Standard Form but with $\lambda = 3\lambda^*$ \subref{fig:realWLS}
Weighted Least Squares \subref{fig:realDIR} reconstruction with
Homogeneous Dirichlet boundary conditions.} \label{fig:realRecPS}
\end{figure*}

There are several test data sets in the literature for Photometric
Stereo.  They are usually photos taken of plaster casts, which are
typically very good approximations to Lambertian surfaces.  As a
consequence, the Photometric Stereo technique yields good
approximations to the gradient of the surface.  Figure
\ref{fig:mozart} shows the reconstruction results of the so-called
``Mozart" data-set.  The results presented in~\cite{agrawal2006}
for the same date set demonstrate that until now this data set
could have been considered challenging.  However, with the global
least squares approach presented here this data set can be
considered almost trivial due to the Lambertian nature of the test
surface. \\
What is, however, far more difficult to reconstruct a surface
which is non-Lambertian with various different surface textures.
Figure \ref{fig:measuredSurf} shows such a surface whose gradient
has been measured via Photometric Stereo.  The surface is metallic
and has irregular textures such as rust.  Thus, the Lambertian
assumption will lead to systematic errors in the gradient
computation.  It is precisely for these such cases that
regularization techniques such as Tikhonov Regularization have
been developed -- for so-called ill-posed problems.
Figure \ref{fig:realRecPS} shows the reconstructions of the
surface using the methods proposed in this paper.  Figure
\ref{fig:realGLS} shows the GLS reconstruction, which exhibits a
global bending due to inhomogeneities in lighting, among other
sources of bias.  In Figure \ref{fig:realSGBP}, spectral
reconstruction with a band-pass filter is used to remove this
bending effect; specifically, all bi-cubic terms are removed from
the reconstruction to remove the bending, and the high-frequency
components are removed to suppress Gaussian-like noise. Similarly,
Tikhonov regularization can be used for the same purpose; in
Figure \ref{fig:realTikLam}, Tikhonov regularization in standard
form has been used, whilst selecting the regularization parameter,
$\lambda_1$, using the L-curve method. To show the difficulty in
determining the regularization parameter, the result using
$\lambda_2=3\lambda_1$ is shown in Figure \ref{fig:realTik3Lam}.
Clearly this simple change produces a much stronger flattening
effect, and the reconstruction appears to be much better in
comparison to the original surface. The effects of the weighted
least squares solution are more subtle; using an inverted
Gaussian-bell like weighting function, the reconstruction errors
in the centre of the image are considered more critical. The
weighted solution in this case is shown in Figure
\ref{fig:realWLS}, which exhibits systematic differences to the
GLS reconstruction, and are more pronounced at the boundaries.
Finally, the advantages of applying homogeneous Dirichlet boundary
conditions to the reconstruction are shown in Figure
\ref{fig:realDIR}; with the edges of the reconstructed surface
essentially simply supported, the low-frequency bending of the
surface is almost completely suppressed.  Such a reconstruction
can be highly effective if the end goal of reconstruction is to
automatically read the code stamped on the steel block.
\section{Conclusion}
This paper presented a framework based on the Sylvester Equation
for direct surface reconstruction methods from gradient fields
with state-of-the-art forms of regularization.  The new algorithms
are several orders of magnitudes faster than previous methods due
to the efficient solution of Sylvester Equations.  A trivial
extension of the Framework, would be to combine the various forms
of regularization, i.e., Spectral Methods combined with Tikhonov
regularization (this has been omitted in the interest of
conciseness). Future work will be to further accelerate the
solution of the Sylvester Equations; clearly, these too are
largely structured and/or sparse.  It should also be noted that
the Sylvester Equations presented here can be partially solved
off-line, and hence may lead to real-time implementations.  The
methods presented here represent the first viable methods for
real-time Photometric Stereo, where regularization is essential,
such as in any Industrial Applications.
\appendix
%
%
\section{Differentiation of a Frobenius Norm w.r.t. a Matrix}
\label{app:diffFrob}
Used frequently throughout this paper is the derivative of the
squared Frobenius norm of the general form,
\begin{equation}
f\left(\M{X}\right) = \left\| \M{A} \M{X} \M{B} - \M{C}
\right\|_{\textrm{F}}^2,
\end{equation}
with respect to the matrix $\M{X}$.
To obtain a formula for the derivative, we firstly define the
derivative of the scalar valued function $f=f\left(\M{X}\right)$
with respect to the $m\times n$ matrix $\M{X}$ as the
%
 matrix of
partial derivatives~\cite{schoenemann1985},
\begin{equation}
\pdiff{f}{\M{X}} = \left[\pdiff{f}{x_{ij}}\right].
\end{equation}
 That is, an $m\times n$ matrix whose $i$-$j$ entry is
the partial derivative of $f\left(\M{X}\right)$ with respect to
the entries $x_{ij}$ of $\M{X}$. The derivative of the Frobenius
norm with respect to the matrix $\M{X}$ is obtained by using the
matrix trace definition of the Frobenius norm, i.e.,
\begin{equation}
f\left(\M{X}\right) = \tr{ \left(\M{A} \M{X} \M{B} -
\M{C}\right)\left(\M{A} \M{X} \M{B} - \M{C}\right)^{\textrm{T}} }
\end{equation}
whereupon expanding yields,
\begin{equation}
f\left(\M{X}\right) = \tr{ \M{A} \M{X} \M{B} \MT{B} \MT{X} \MT{A}
 - \M{C} \MT{B} \MT{X} \MT{A} - \M{A} \M{X} \M{B} \MT{C} + \M{C}\MT{C}
 }.
\end{equation}
%
%
Thus, noting that
\begin{equation}
\pdiff{\M{X}}{x_{ij}} = \V{e}_i\VT{e}_j,
\end{equation}
the derivative of the function with respect to the entry $x_{ij}$
is,
\begin{eqnarray}
\pdiff{f}{x_{ij}} & = & \mathop{\mathrm{trace}}\left(
 \M{A} \V{e}_i\VT{e}_j \M{B} \MT{B} \MT{X} \MT{A} + \M{A} \M{X} \M{B} \MT{B} \V{e}_j\VT{e}_i \MT{A} \right. \nonumber \\
& - & \left.
  \M{C} \MT{B} \V{e}_j\VT{e}_i \MT{A} - \M{A} \V{e}_i\VT{e}_j \M{B} \MT{C}
  \right).
\end{eqnarray}
Due to the definition of the trace, and the symmetry of its
argument, we have,
\begin{equation}
\pdiff{f}{x_{ij}} = 2\tr{ \M{A} \V{e}_i\VT{e}_j \M{B} \MT{B}
\MT{X} \MT{A} - \M{A} \V{e}_i\VT{e}_j \M{B} \MT{C} }
\end{equation}
or more simply,
\begin{equation}
\pdiff{f}{x_{ij}} = 2\tr{ \M{A} \V{e}_i\VT{e}_j \MT{M} },
\end{equation}
where the matrix $\M{M}$ is the placeholder,
\begin{equation}
\M{M} = \M{A}\M{X}\M{B}\MT{B} - \M{C}\MT{B}.
\end{equation}
\clearpage
\noindent Thus by appropriately indexing the matrices $\M{A}$ and
$\M{M}$ we have,
\begin{equation}
\pdiff{f}{x_{ij}} = 2 \sum_{k} a_{ki} m_{kj}.
\end{equation}
The matrix whose $i$-$j$ entry is this expression is obtained from
the rules of matrix multiplication, and hence,
\begin{equation}
\pdiff{f}{\M{X}} = 2 \MT{A} \M{M}.
\end{equation}
By replacing $\M{M}$ we obtain the desired formula for the
derivative of the Frobenius norm,
\begin{equation}
\pdiff{}{\M{X}} \left\| \M{A} \M{X} \M{B} - \M{C}
\right\|_{\textrm{F}}^2
 = 2\MT{A} \left( \M{A} \M{X} \M{B} - \M{C} \right)\MT{B}
\end{equation}
This identity can be derived using the methods developed in
Sch\"{o}nemann~\cite{schoenemann1985}. All derivatives in this
paper can be found using this identity with special cases such as
$\M{A}=\M{I}$ or $\M{B}=\M{I}$.
%
%
\section*{Acknowledgment}
The authors would like to thank Georg Jaindl for acquiring the
Photometric Stereo images~\cite{jaindl2009}.
%
%
\bibliographystyle{spmpsci}      
\bibliography{BIB/grad2surf-Converted}   
\end{document}